\documentclass[preprint]{elsarticle}
\usepackage{lineno,hyperref}
\modulolinenumbers[5]
\journal{}
\usepackage{setspace}

\usepackage{geometry}

\newgeometry{vmargin={25mm}, hmargin={20mm,20mm}}
\usepackage{color}
\usepackage{amsmath}
\usepackage{subcaption}
\usepackage{amsfonts}

\graphicspath{{Figs/}}
\begin{document}

\begin{frontmatter}

\title{Energy-recurrence Breakdown and Chaos in Disordered Fermi-Pasta-Ulam-Tsingou Lattices}

%
%
%

	\author[au1,au2]{Zulkarnain}
\author[au3]{H.\ Susanto\corref{cor1}}
\ead{hadi.susanto@ku.ac.ae}
\author[au1]{C.\ G.\ Antonopoulos}

\address[au1]{Department of Mathematical Sciences, University of Essex, Wivenhoe Park, Colchester CO4 3SQ, United Kingdom}
\address[au2]{Department of Mathematics, Faculty of Mathematics and Natural Sciences, Universitas Riau, Pekanbaru 28293, Indonesia}
\address[au3]{Department of Mathematics, Khalifa University, PO Box 127788, Abu Dhabi, United Arab Emirates}
\cortext[cor1]{Corresponding author}

\begin{abstract}
In this paper, we consider the classic Fermi-Pasta-Ulam-Tsingou system as a model of interacting particles connected by harmonic springs with a quadratic nonlinear term (first system) and a set of second-order ordinary differential equations with variability (second system) that resembles Hamilton's equations of motion of the Fermi-Pasta-Ulam-Tsingou system. In the absence of variability, the second system becomes Hamilton's equations of motion of the Fermi-Pasta-Ulam-Tsingou system (first system). Variability is introduced to Hamilton's equations of motion of the Fermi-Pasta-Ulam-Tsingou system to take into account inherent variations (for example, due to manufacturing processes), giving rise to heterogeneity in its parameters. We demonstrate that a percentage of variability smaller than a threshold can break the well-known energy recurrence phenomenon and induce localization in the energy normal-mode space. However, percentage of variability larger than the threshold may make the trajectories of the second system blow up in finite time. Using a multiple-scale expansion, we derive analytically a two normal-mode approximation that explains the mechanism for energy localization and blow up in the second system. We also investigate the chaotic behavior of the two systems as the percentage of variability is increased, utilising the maximum Lyapunov exponent and Smaller Alignment Index. Our analysis shows that when there is almost energy localization in the second system, it is more probable to observe chaos, as the number of particles increases.
\end{abstract}


\end{frontmatter}


\section{Introduction}

Debye suggested that thermal conductivity in a crystal is a consequence of atom vibrations in the lattice \cite{peierls1929kinetischen,debye1912theorie}. To model thermalization processes in physical media, Fermi, Pasta, Ulam, with Tsingou's help running the computer simulations \cite{fermi1955alamos}, considered a system of particles connected by harmonic springs with a quadratic nonlinear term, i.e., the so-called FPUT lattice, that is fixed at both ends.

A purely linear dynamics of the springs keeps energy, given to a single normal mode, localized in that mode. However, introducing nonlinear interactions, one would expect that energy introduced to one normal mode, would slowly spread to other normal modes, until the system reaches a state of equipartition of energy, i.e., the system relaxes to a thermal equilibrium. 
Contrary to this expectation, Fermi, Pasta, Ulam, and Tsingou observed in their seminal paper \cite{fermi1955alamos} in 1955 recurrences of the energy to its initial state, the so-called FPUT recurrences, a phenomenon that led to numerous discoveries in mathematics and physics \cite{ford1992fermi,berman2005fermi,dauxois2005fermi,gallavotti2007fermi,porter2009fermi} thereafter. 

The recurrence phenomenon was explained by Zabusky and Kruskal \cite{zabusky1965interaction} in real space, who derived the integrable Korteweg-de Vries equation from the continuum limit of the FPUT lattice. Introducing energy into one normal mode with wave number $k$ is nothing else but taking the sinusoidal initial condition in real space. As time evolves, the state breaks into a series of localized solutions, i.e., solitons that move and interact with the fixed ends, i.e. boundaries. Upon interacting with the fixed ends, the solitons bounce back and return to their initial positions, i.e., giving rise to FPUT recurrences. Another explanation to the inefficient energy transfer among normal modes was provided by Izrailev and Chirikov \cite{izrailev1966statistical} who used the concept of the overlap of nonlinear resonances. They associated equipartition of energy with dynamical chaos and were able to estimate a threshold that separates regular from chaotic dynamics.

Another direction in the study of the FPUT lattice is that of heat conductivity in the presence of disorder. The main interest is in its interplay with nonlinearity. For harmonic disordered systems, all eigenmodes of the infinite system, i.e., Anderson modes, are known to be localized and form a complete basis \cite{anderson1958absence}. As a linear superposition of Anderson modes, an initially localized wave in the infinite chain will remain localized at any time. Whether this behavior changes qualitatively by the introduction of nonlinearity is still an open question (see for example \cite{li2001can,dhar2008heat,zhu2021effects} and references therein). Disorder can be introduced in the form of uniformly distributed random variation of particle masses \cite{payton1967lattice,li2001can}, linear coupling constants between nearest neighbours \cite{lepri2010asymptotic}, or in the nonlinearity coefficients \cite{pikovsky2015first}. Recently, by viewing FPUT lattices as systems of masses coupled with nonlinear springs, Nelson et al.\ \cite{porter2018} incorporated heterogeneity on a one-dimensional FPUT array to take into account uncertainties (i.e., in the masses, the spring constants, or the nonlinear coefficients) during the manufacturing process of such physical systems. They demonstrated numerically that tolerances degrade the observance of recurrences, often leading to a complete loss in moderately-sized arrays. Such a variability may therefore provide a plausible explanation to little experimental evidence on FPUT energy recurrences. 

Here, we consider the problem of heterogeneous FPUT systems studied in \cite{porter2018}. 
In our work, we perform numerical simulations in great detail to understand the breakdown of FPUT recurrences in the model. Indeed, we observe recurrence degradation, where the energy peak of the lowest normal mode is decreasing subsequently. For percentage of variability smaller than a threshold that we derive, the energy is then localized in the few lowest normal modes. The authors in \cite{flach2005q,flach2006q} considered non-equipartition of energy among normal modes and studied time-periodic states that are exponentially localized in the $k$- (or $q$- in \cite{flach2005q,flach2006q}) space of normal modes. Such time-periodic states are referred to as $q$-breathers. Variability in FPUT lattices therefore leads to $q$-breathers. In our work, by transforming the FPUT system into another system in the normal-mode space and considering a two normal-mode approximation, we provide a qualitative explanation for the disappearance of FPUT recurrences. In this approximation, $q$-breathers are periodic solutions centered around an equilibrium point (i.e., time-independent solutions) in the $q$-space of normal modes.

We also perform long-term numerical integrations and compute the maximum Lyapunov exponent (mLE) \cite{benettin1980lyapunov} and Smaller Alignment Index (SALI) \cite{skokos2003does,skokos2004detecting} to show that the trajectories of the heterogeneous system become quickly chaotic, as the number of particles in the system increases for the same percentage of variability. In homogeneous FPUT lattices (i.e., in the absence of variability), when recurrences occur, the system reaches a metastable state \cite{ponno2005korteweg,bambusi2006metastability}, where only few (low $k$) normal modes share the total energy of the system. However, it has also been shown that a rather weak diffusion takes place in the highest normal modes of the spectrum \cite{ponno2011two} that gradually leads to equipartition of energy \cite{berchialla2004exponentially}. Using the result in \cite{izrailev1966statistical}, this weak diffusion process implies weak chaos \cite{antonopoulos2011weak}. Our work shows that variability enhances the chaotic dynamics of the system.

In this work, we also show that for percentages of variability bigger than a threshold, solutions may blow up in finite time. Using the same two normal-mode approximation, we have been able to explain the blow up phenomenon. A bifurcation analysis is further provided that yields a variability threshold for the blow up of solutions.

The paper is organised as follows: In Sec.\ \ref{sec2}, we review the original FPUT lattice with a quadratic nonlinearity (i.e., the FPUT$-\alpha$ system) and discuss energy recurrence. We introduce the governing equations of motion in the presence of parameter variability in Sec.\ \ref{sec3}. The phenomena of recurrence breakdown and blow up of solutions are reported in the same section. In Sec.\ \ref{sec4}, a two normal-mode approximation in the normal mode space is derived using multiple-scale analysis. Our analytical results explain why energy recurrences breakdown when variability is introduced and provide a qualitative reason why solutions blow up in finite-time after a variability threshold. In Sec.\ \ref{sec5}, we discuss chaos in the FPUT$-\alpha$ system with or without variability and the mLE and SALI methods that we use to discriminate between ordered and chaotic trajectories. Finally, we conclude our study and discuss future work in Sec.\ \ref{conc}.

\section{Mathematical model and dynamics of FPUT$-\alpha$ lattices}\label{sec2}

The Hamiltonian of the FPUT$-\alpha$ system is given by
\begin{equation}
	H(x,p)=\frac{1}{2}\sum_{j=0}^N p_j^2 + \sum_{j=0}^N \frac{1}{2} \left( x_{j+1}-x_j \right)^2 + \frac{\alpha}{3} \left( x_{j+1}-x_j \right)^3=E,
	\label{Hamilt}
\end{equation}
where fixed boundary conditions $x_0 = x_{N+1}=0$ and $p_0=0$ are considered. In this context, $\alpha\geq0$ is the nonlinear coupling strength and $E$ the total, fixed, energy of the system. By viewing the FPUT$-\alpha$ lattice as a model of particles coupled with springs, $x_j(t)$ represents the relative displacement of the $j$th-particle from its equilibrium position at any time $t$ and $p_j(t)$ its corresponding conjugate momentum at any time $t$. The equations of motion that result from Hamiltonian \eqref{Hamilt} (i.e., Hamilton's equations of motion) are then given by
\begin{align}
	\ddot{x}_j= & (x_{j+1}-x_j)+\alpha(x_{j+1}-x_j)^2-(x_{j}-x_{j-1}) - \alpha(x_{j}-x_{j-1})^2.
	\label{eqnmotions}
\end{align}

Working in the real space $x_j$ and $p_j$, one can express Eqs. \eqref{eqnmotions} in the normal-mode space $Q_j$ and $P_j$. This can be done by writing the position $x_j(t)$ as a superposition of eigenvectors of the linear equation. Using the normal mode transformation, 
\begin{align}
	\mathbf{x} &=A\mathbf{Q},\quad \mathbf{p} =A\mathbf{P}, \label{Pk}
\end{align}
where $\mathbf{x}=[x_1~x_2~\ldots~x_N]^T$, $\mathbf{p}=[p_1~p_2~\ldots~p_N]^T$, $\mathbf{Q}=[Q_1~Q_2~\ldots~Q_N]^T$, $\mathbf{P}=[P_1~P_2~\ldots~P_N]^T$, and
\begin{equation}
	A=\sqrt{\frac{2}{N+1}}\begin{bmatrix}
		\sin \left(\frac{\pi}{N+1} \right ) & \sin \left(\frac{2\pi}{N+1} \right ) & \dots & \sin \left(\frac{N\pi}{N+1} \right )\\ 
		\sin \left(\frac{2\pi}{N+1} \right )& \sin \left(\frac{4\pi}{N+1} \right ) & \dots & \sin \left(\frac{2N\pi}{N+1} \right )\\ 
		\vdots & \vdots & \ddots & \vdots\\ 
		\sin \left(\frac{N\pi}{N+1} \right )& \sin \left(\frac{2N\pi}{N+1} \right ) & \dots & \sin \left(\frac{N^2\pi}{N+1} \right )
	\end{bmatrix},
	\label{A_matrix}
\end{equation}
%
the Hamiltonian \eqref{Hamilt} becomes
\begin{align*}
	H&=\frac{1}{2}\sum_{k=1}^N \left(P_k^2+\omega_k^2Q_k^2 \right)+\alpha H_3(Q_1,Q_2,\ldots,Q_N),
\end{align*}
for some nonlinear function $H_3$, where 
\begin{equation}
	\omega_k=2\sin\left(\frac{k\pi}{2(N+1)}\right).
	\label{wk}
\end{equation}
In this framework, $\mathbf{Q}$ represents the amplitude of the normal mode, while $\mathbf{P}$ its velocity. The energy of normal mode $k$ for $\alpha=0$ can then be defined by 
\begin{align}
	E_k&=\frac{1}{2} \left(P_k^2+\omega_k^2Q_k^2 \right).\label{nmE}
\end{align}

Substituting Eq.\ \eqref{Pk} into Eq.\ \eqref{eqnmotions}, we obtain the equations of motion in normal-mode coordinates as
\begin{align}
	\ddot{\mathbf{Q}} &=D\mathbf{Q}+A^{-1}\mathbf{F(Q)},
	\label{eqnmotionsnmcoord}
\end{align} 
where 
\begin{align*}
	D=\begin{bmatrix}
		-\omega_1^2 & 0 & \ldots & 0\\ 
		0& -\omega_2^2 &  & 0\\ 
		\vdots&  & \ddots & \vdots\\ 
		0& 0 & \ldots & -\omega_N^2
	\end{bmatrix},\;\;
	\mathbf{F(Q)} =\begin{bmatrix}
		f_1(\mathbf{Q})\\ 
		f_2(\mathbf{Q})\\ 
		\vdots\\
		f_N(\mathbf{Q})
	\end{bmatrix}
\end{align*}
and $A^{-1}$ is the inverse matrix of $A$, given by Eq.\ \eqref{A_matrix}.

In their seminal work, Fermi, Pasta, Ulam and Tsingou \cite{fermi1955alamos} excited the lowest possible normal mode, i.e., the mode with $k=1$. The initial conditions of Eqs. \eqref{eqnmotions} are then
\begin{align}
	p_j&=0, \quad x_j=\sin \left( \frac{\pi i}{N+1}\right),\;j=1,2,\ldots,N, \label{initxi}
\end{align} 
which are equivalent to solving Eq.\ \eqref{eqnmotionsnmcoord} with $Q_1=\sqrt{(N+1)/2}$, $Q_k=0$ for $k=2,3,\dots,N,$ and $\dot{Q}_k=0$ for $k=1,2,\dots,N$.

\begin{figure}[h]
	\begin{subfigure}{.5\textwidth}
		\centering
		\includegraphics[width=0.9\linewidth]{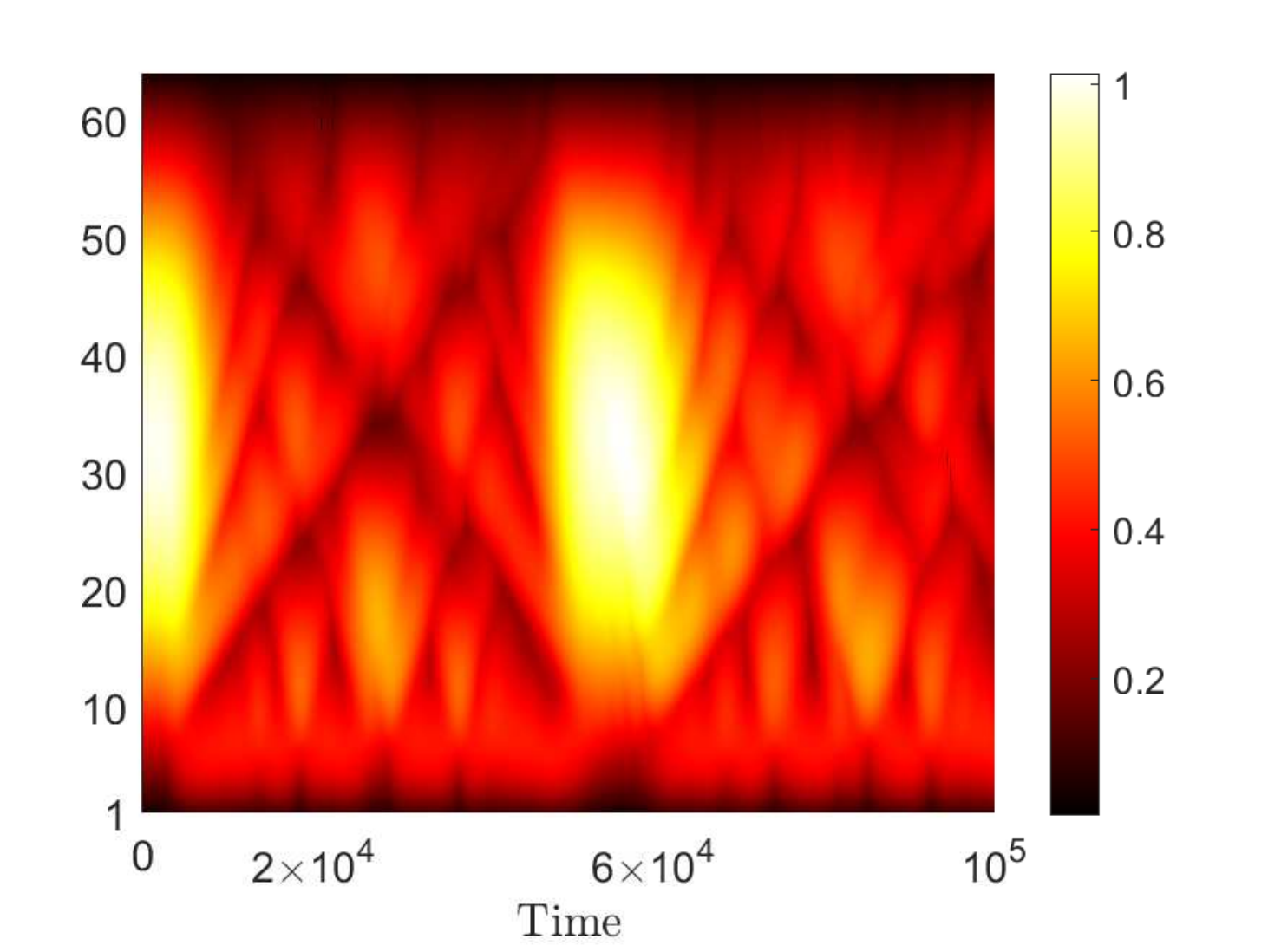}
		\caption{}
		\label{}
	\end{subfigure}%
	\begin{subfigure}{.5\textwidth}
		\centering
		\includegraphics[width=0.9\linewidth]{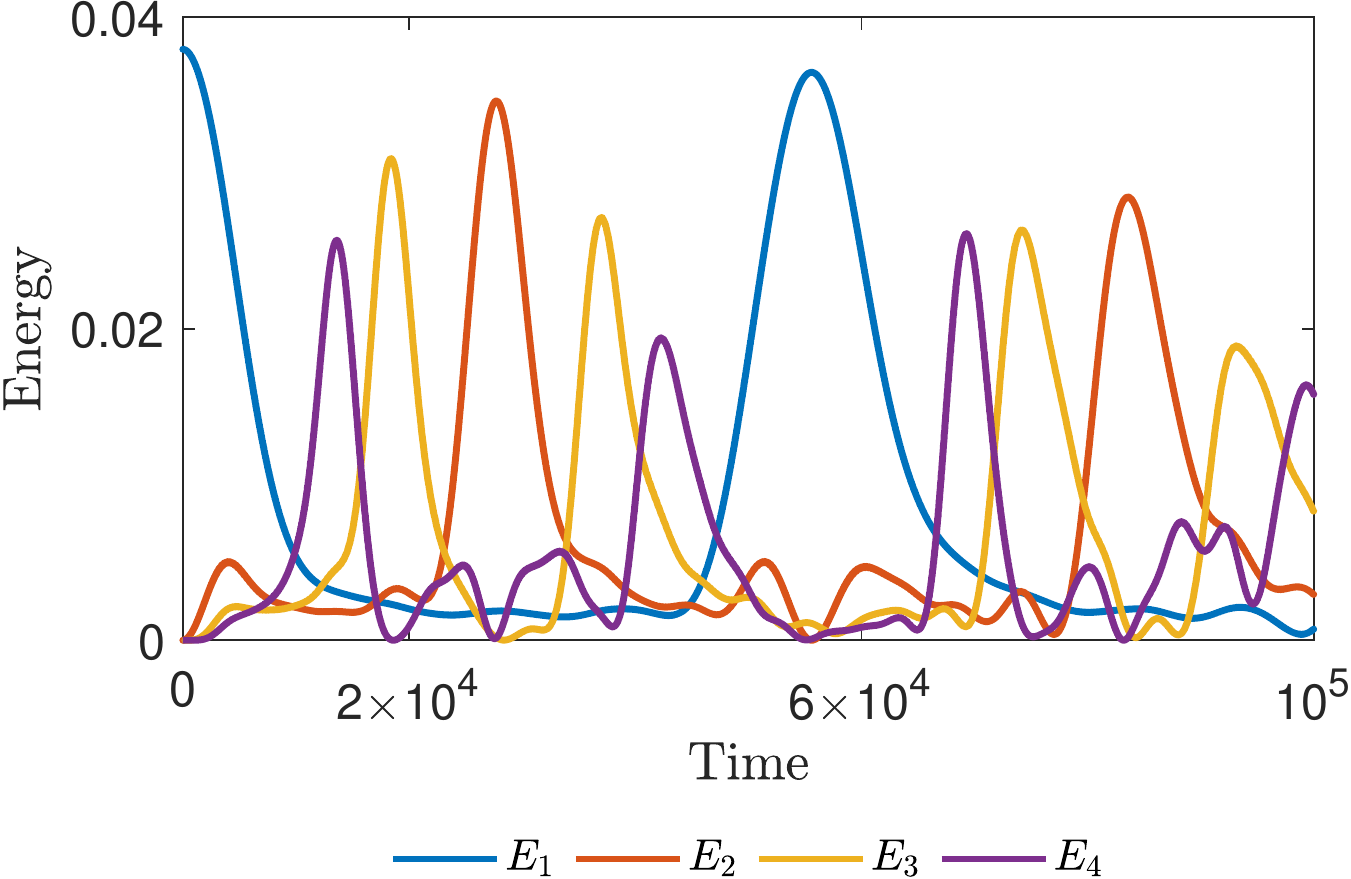}
		\caption{}
		\label{}
	\end{subfigure}
	\caption{FPUT recurrences of the Hamiltonian system \eqref{Hamilt}. (a) Dynamics of $x_j(t)$ using the initial condition in \eqref{initxi}, where $E=0.03795$. Panel (a) shows the top view of the oscillation envelope of $x_j(t)$ in time. (b) Energy of the first four normal modes in the dynamics shown in panel (a). Note in panel (b) how almost all of the energy returns to the first normal mode at around $t=6\times10^4$, i.e., the appearance of an FPUT recurrence. Here we have used $N=64$ in the computations in both panels. The range of values in the vertical axis in panel (a) is between 1 and $N=64$.}\label{fig:fig1}
\end{figure}

In Fig.\ \ref{fig:fig1}, we plot the dynamics of $x_j$ of system \eqref{Hamilt} for $N=64$, for the initial condition \eqref{initxi}, where $E=0.03795$. Panels (a) and (b) show the dynamics of $x_j(t)$ in real space (infact, what is shown is the oscillation envelope) and the normal mode energy of the first four normal modes of the FPUT$-\alpha$ system \eqref{Hamilt}, respectively.

In their seminal paper  \cite{fermi1955alamos}, Fermi, Pasta, Ulam and Tsingou expected that the energy $E$, which was initially used to excite the lowest normal mode only (i.e., $k=1$), would slowly drift to the other normal modes until the system reaches thermalization, as predicted by Statistical Mechanics. Surprisingly, the numerical experiment showed that that was not the case and that after several periods of the evolution of the mode, almost all energy in the system returned to the first normal mode that was excited initially. The authors witnessed the so-called FPUT recurrences. An example of such recurrences is given in Fig.\ \ref{fig:fig1} for $N=64$ and $E=0.03795$.

\section{Disordered FPUT lattices}\label{sec3}
The authors in \cite{porter2018} proposed various disordered FPUT$-\alpha$ systems that include tolerances $v$ into each particle as a result of variability in a manufacturing process. In their study, they claimed to incorporate tolerances into the system in different ways based on manufacturing constraints. They first introduced the following system with heterogeneity 
\begin{align}
	\ddot{x}_j= & (v_{j+1}x_{j+1}-v_jx_j)+\alpha(v_{j+1}x_{j+1}-v_jx_j)^2-(v_jx_{j}-v_{j-1}x_{j-1}) - \alpha(v_jx_{j}-v_{j-1}x_{j-1})^2,\label{dis1}
\end{align}
where, again $\alpha\geq0$ is the nonlinear coupling strength. It can be shown that this system admits the Hamiltonian function 
\begin{equation}
	H(x,p)=\frac{1}{2}\sum_{j=0}^N \frac{p_j^2}{v_j} + \sum_{j=0}^N \frac{1}{2} \left( v_{j+1}x_{j+1}-v_jx_j \right)^2 + \frac{\alpha}{3} \left( v_{j+1}x_{j+1}-v_jx_j \right)^3=E.
	\label{Hamilt2mod}
\end{equation}
The disorder is inserted in a symmetric way between the linear and nonlinear coupling. Particularly, the variabilities $v_j$ were generated randomly from a Gaussian distribution, that is for a tolerance $\tau\%$, the values of $v_j$ were drawn from a Gaussian distribution with mean 1 and standard deviation $\sigma=1/3 \times 0.01\tau$. Therefore, the values of $v_j$ would lie in the interval $[1-0.01\tau,1+0.01\tau]$.

The authors in \cite{porter2018} considered also the case where the variabilities $v_j$ are present only in the nonlinear coupling terms, resulting in the following system of second-order ordinary differential equations
\begin{align}
	\ddot{x}_j= & (x_{j+1}-x_j)+\alpha(v_{j+1}x_{j+1}-v_jx_j)^2-(x_{j}-x_{j-1}) - \alpha(v_jx_{j}-v_{j-1}x_{j-1})^2.\label{varnonlin}
\end{align}
In this case, the system is no longer Hamiltonian. Then, they showed numerically that incorporating variability in the nonlinear coupling terms only has, for a fixed amount of variability, a comparable effect to incorporating it in only the linear coupling terms. Although in both setups recurrences such as those in Fig.\ \ref{fig:fig1} disappear for large enough tolerance and the energy localizes in the first few normal modes, more energy is transferred to the lower modes in the latter case than in the former.

In this work, we consider the effect of disorder in the second scenario of Eq.\ \eqref{varnonlin}, which is a toy-dyamical system that does not necessarily relate to a real physical system. We have decided to study it as we show in Sec. \ref{sec4}, we can derive a mathematical theory to understand the effect of variability in the localization of energy and its dynamics.

\begin{figure}[h]
	\begin{subfigure}{.5\textwidth}
		\centering
		\includegraphics[ width=0.9\linewidth]{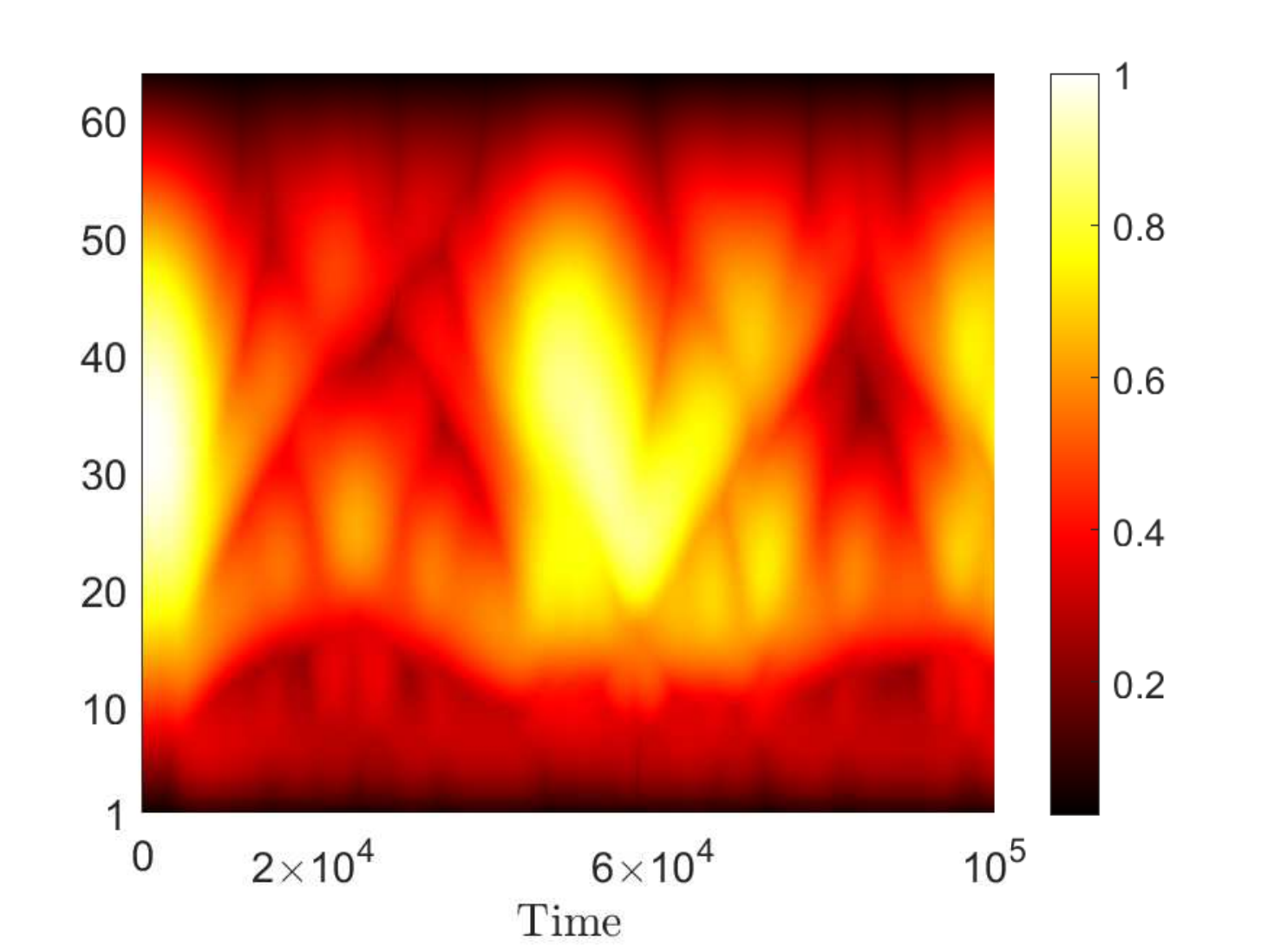}
		\caption{}
		\label{}
	\end{subfigure}%
	\begin{subfigure}{.5\textwidth}
		\centering
		\includegraphics[width=0.9\linewidth]{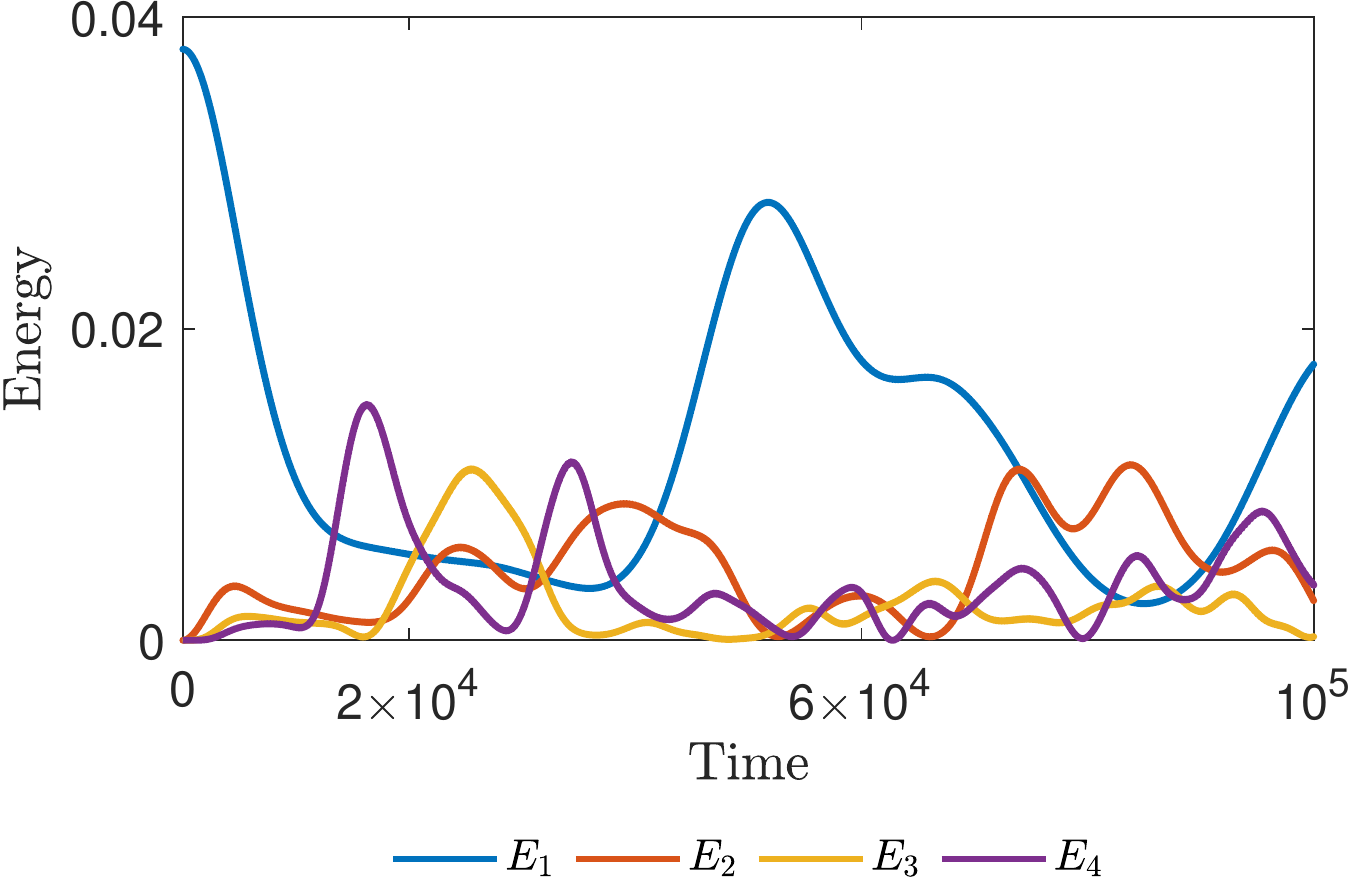}
		\caption{}
		\label{}
	\end{subfigure} \\
	\begin{subfigure}{.5\textwidth}
		\centering
		\includegraphics[width=0.9\linewidth]{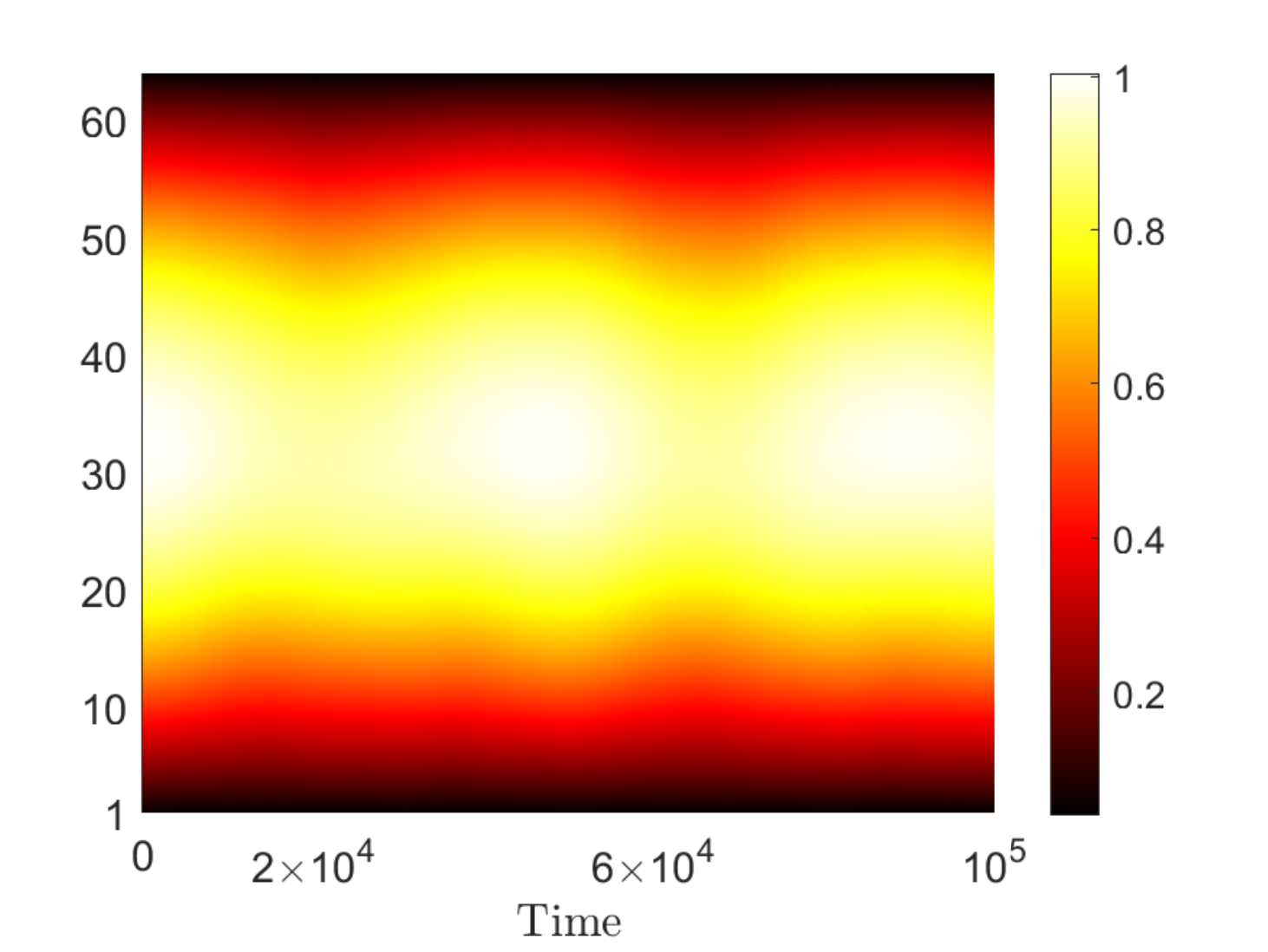}
		\caption{}
		\label{}
	\end{subfigure}%
	\begin{subfigure}{.5\textwidth}
		\centering
		\includegraphics[width=0.9\linewidth]{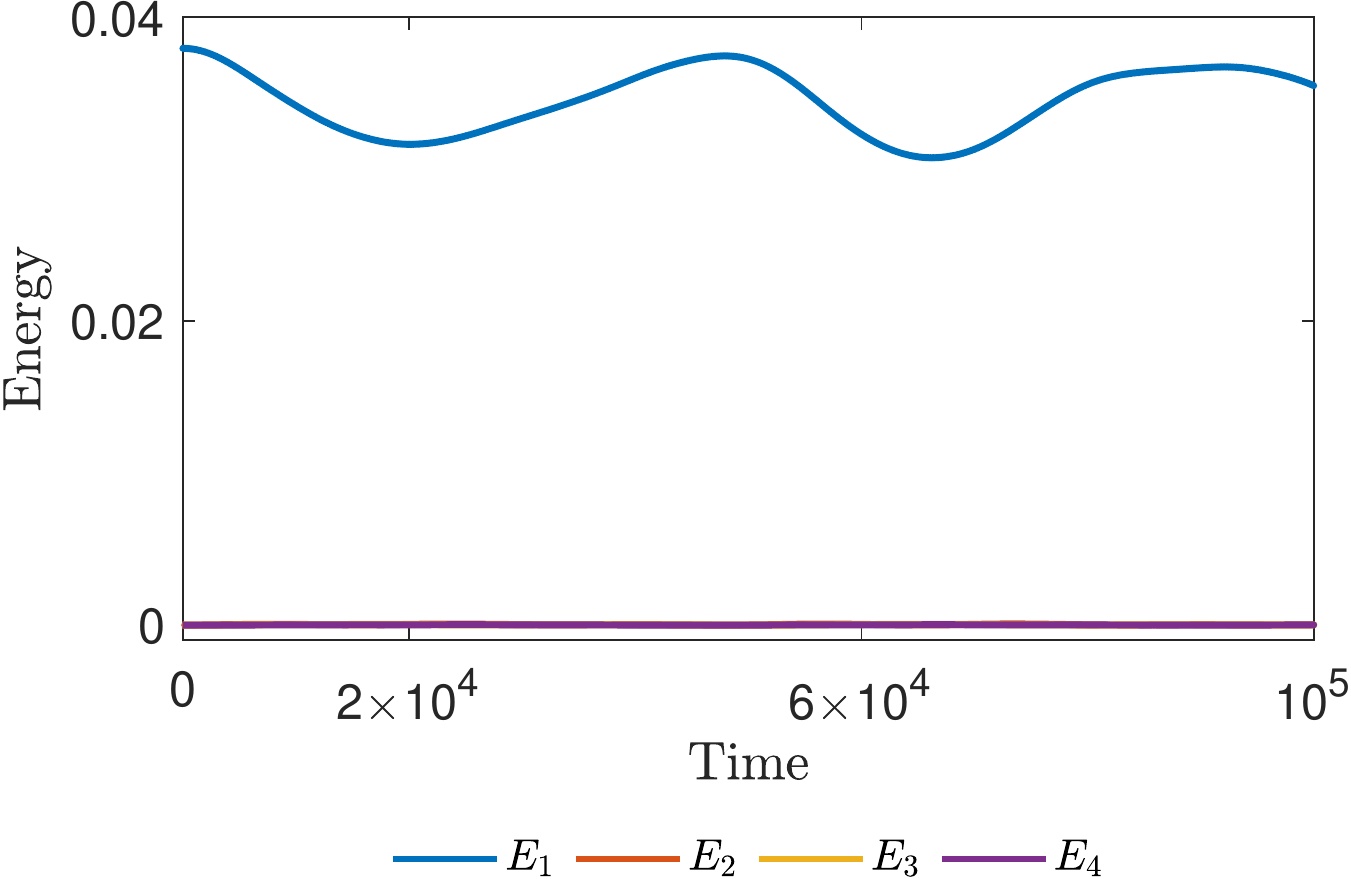}
		\caption{}
		\label{}
	\end{subfigure}
	\caption{Energy recurrences in time for the system in Eq.\ \eqref{varnonlin} for $N=64$, similarly to Fig.\ \ref{fig:fig1} for the Hamiltonian system \eqref{Hamilt}. Panels (a) and (b) are for $\tau=5\%$ tolerance and panels (c) and (d) for $\tau=10\%$. Note that the ranges in the vertical axes in panels (a) and (c) are from 1 to $N=64$.}\label{fig:fig2}
\end{figure}

Throughout the paper, we consider disorder that is generated using the same setup as in \cite{porter2018}. We show in Fig.\ \ref{fig:fig2} the dynamics of $x_j$ and normal mode energies $E_k$ of the first four modes of the system  of equations \eqref{varnonlin} (i.e., for $k=1,\ldots,4$) for $N=64$ particles and two different percentages of tolerance, i.e., for $\tau=5\%$ and $10\%$.

Comparing panels (a) in Figs. \ref{fig:fig1} and \ref{fig:fig2}, one can see that the variability reduces the effectiveness of recurrence, where a subsequent peak of the mode energy $E_1$ is lower than the preceding ones. In \cite{porter2018}, it was reported that for larger variability, the energy transfer from the lowest to the higher ones becomes ineffective, which creates a non-recurrent state, shown in panel (d) in Fig.\ \ref{fig:fig2}. This state is localized in the normal mode space, i.e., it is a $q$-breather \cite{flach2005q,flach2006q}. In other words, disorder promotes the occurrence of $q$-breathers.

\begin{figure}[h]
	\centering
	\includegraphics[width=0.5\textwidth]{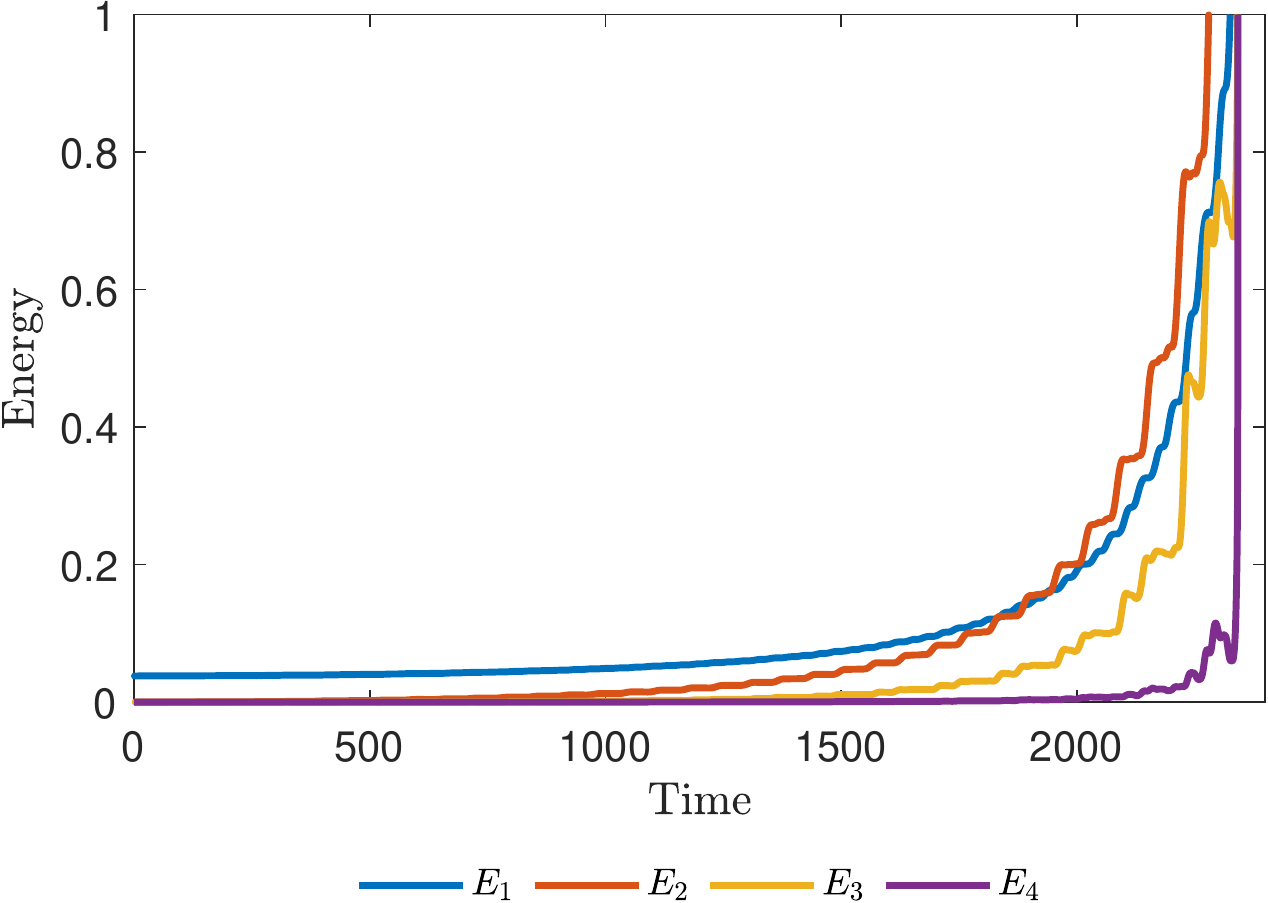}
	\caption{A similar simulation as in panels (b) and (d) in Fig.\ \ref{fig:fig2}, for the system in Eq.\ \eqref{varnonlin} and $N=64$, but for the increased tolerance $\tau=20\%$, where a finite-time blow up of the solution manifests as the abrupt increase of the energy of the first four normal modes around $t=2000$.}\label{fig:fig3}
\end{figure}

In Sec.\ \ref{sec4}, applying a two normal-mode approximation to Eqs. \eqref{eqnmotions} and using multiple-scale expansions, we show that there is a threshold for the percentage of variability $\tau_c\approx10.0749\%$, after which the initial condition \eqref{initxi} may lead to finite-time blow up of the solutions. This is illustrated in Fig.\ \ref{fig:fig3}, where the effect of the blow up is clearly seen in the abrupt increase of the normal-mode energies $E_k$ (see Eq.\ \eqref{nmE}) of the first four normal modes (i.e., for $k=1,\ldots,4$) for $\tau=20\%>\tau_c$.

\begin{figure}[h]
	\centering
	\begin{subfigure}{.3\textwidth}
		\centering
		\includegraphics[width=0.9\linewidth]{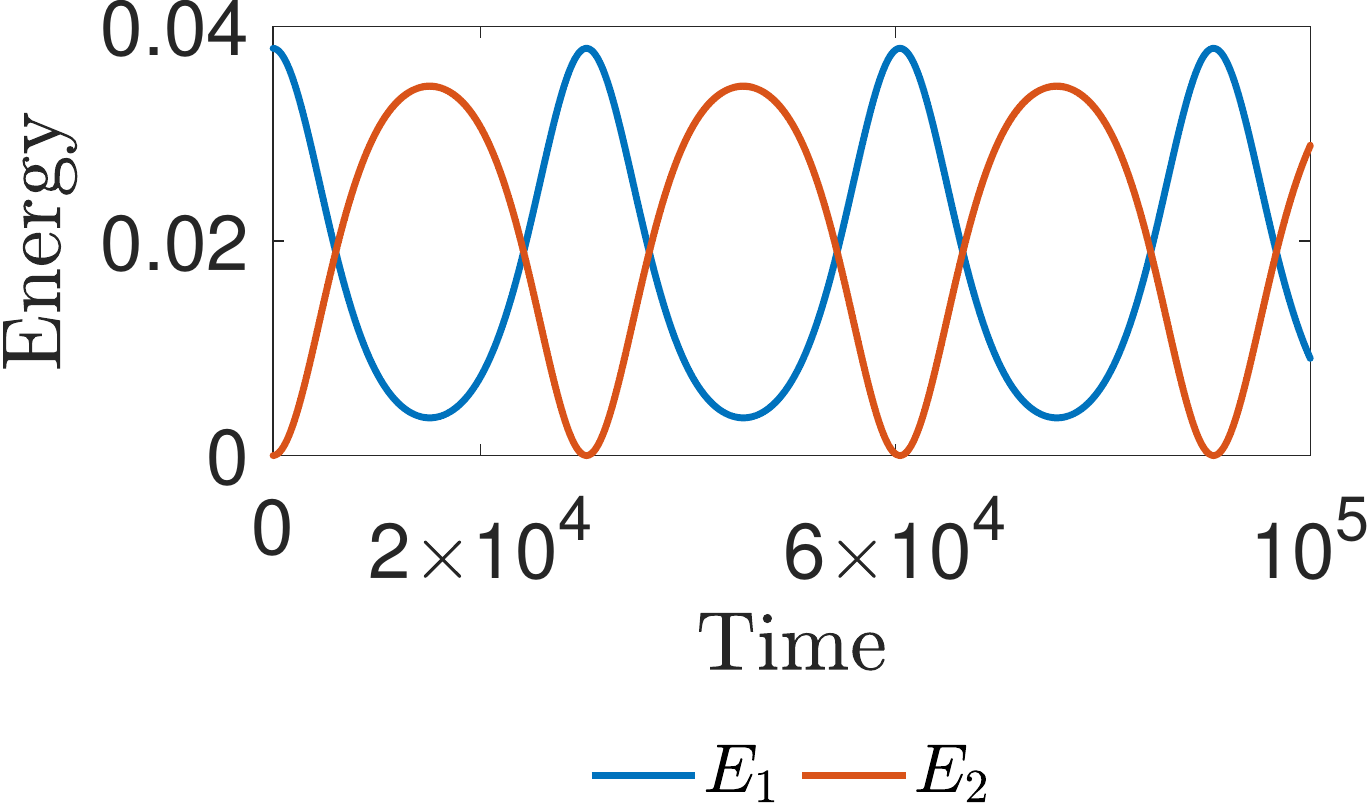}
		\caption{}
		\label{}
	\end{subfigure}%
	\begin{subfigure}{.3\textwidth}
		\centering
		\includegraphics[width=0.9\linewidth]{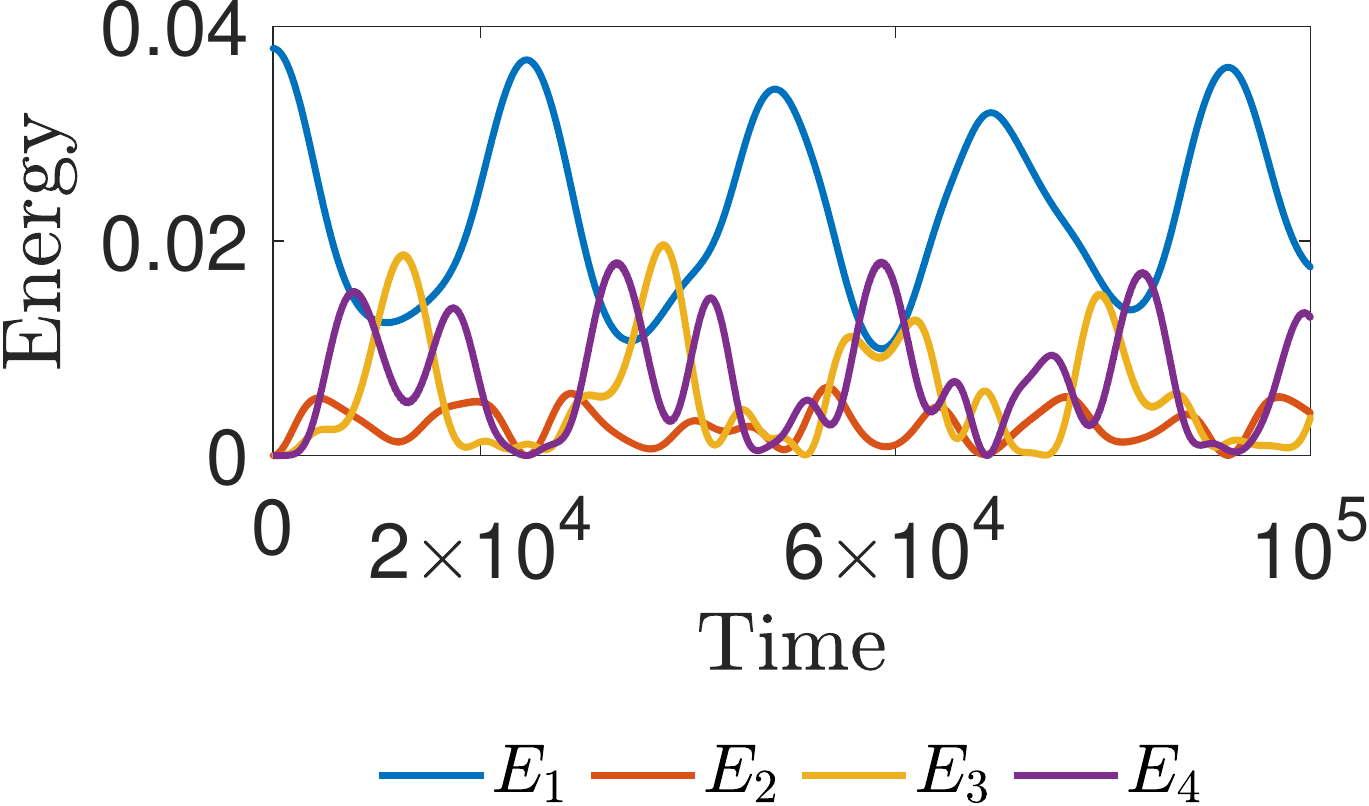}
		\caption{}
		\label{}
	\end{subfigure}%
	\begin{subfigure}{.3\textwidth}
		\centering
		\includegraphics[width=0.9\linewidth]{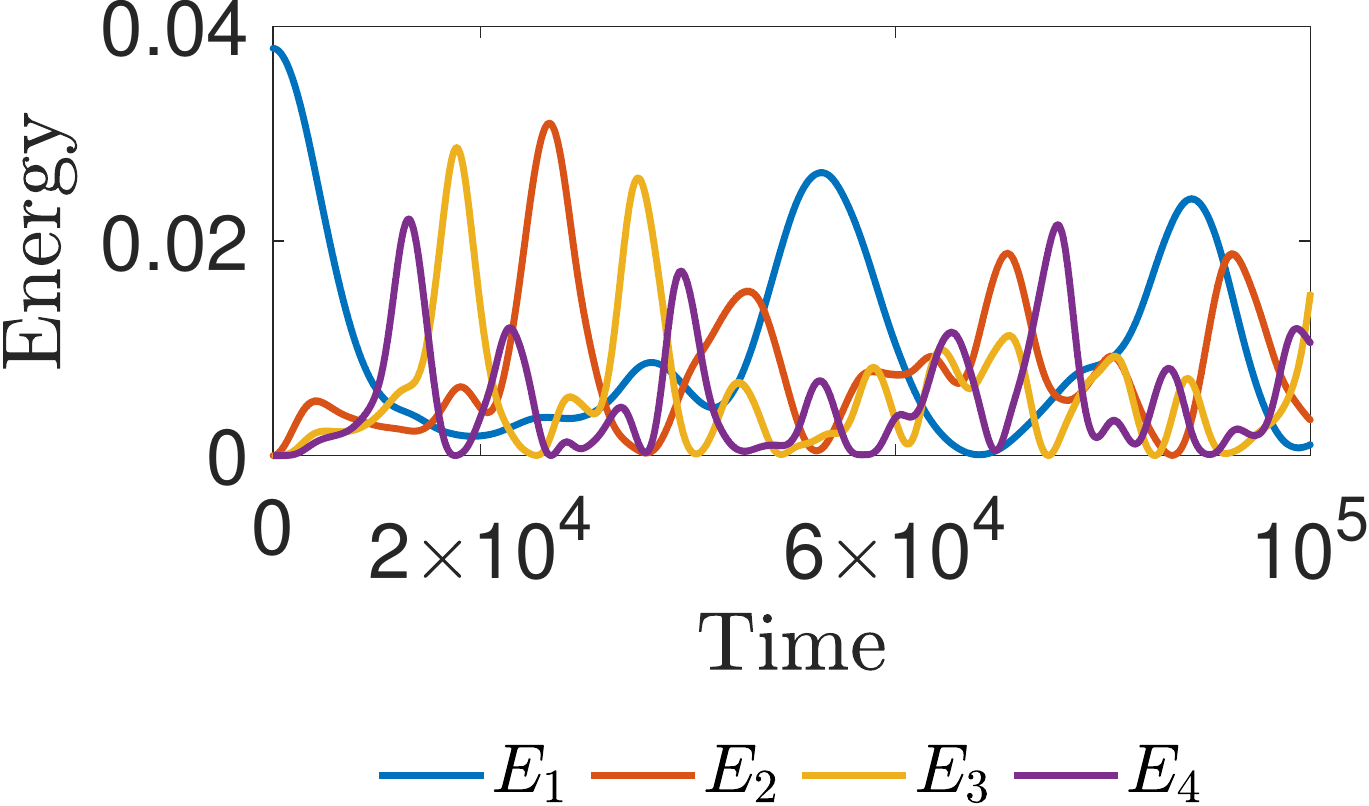}
		\caption{}
		\label{}
	\end{subfigure}\\
	\begin{subfigure}{.3\textwidth}
		\centering
		\includegraphics[width=0.9\linewidth]{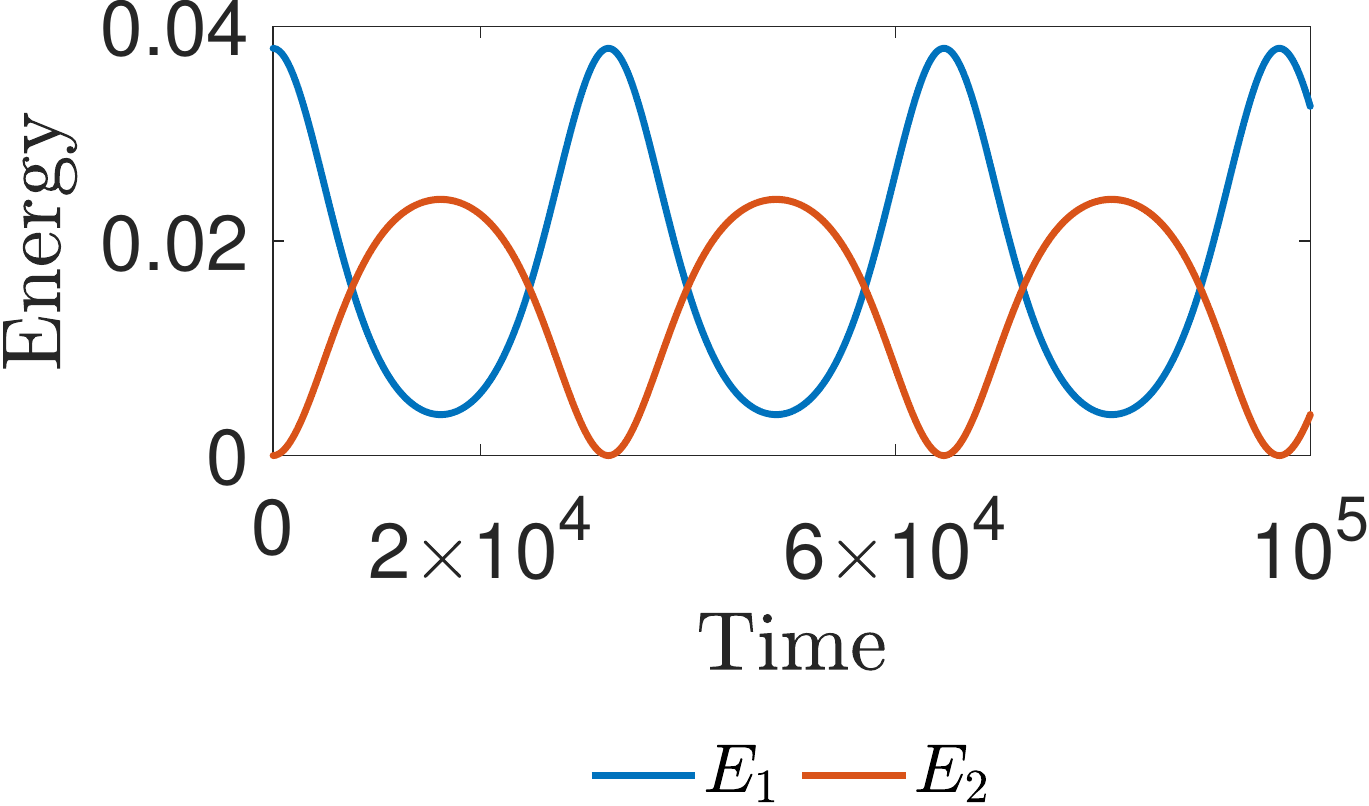}
		\caption{}
		\label{}
	\end{subfigure}%
	\begin{subfigure}{.3\textwidth}
		\centering
		\includegraphics[width=0.9\linewidth]{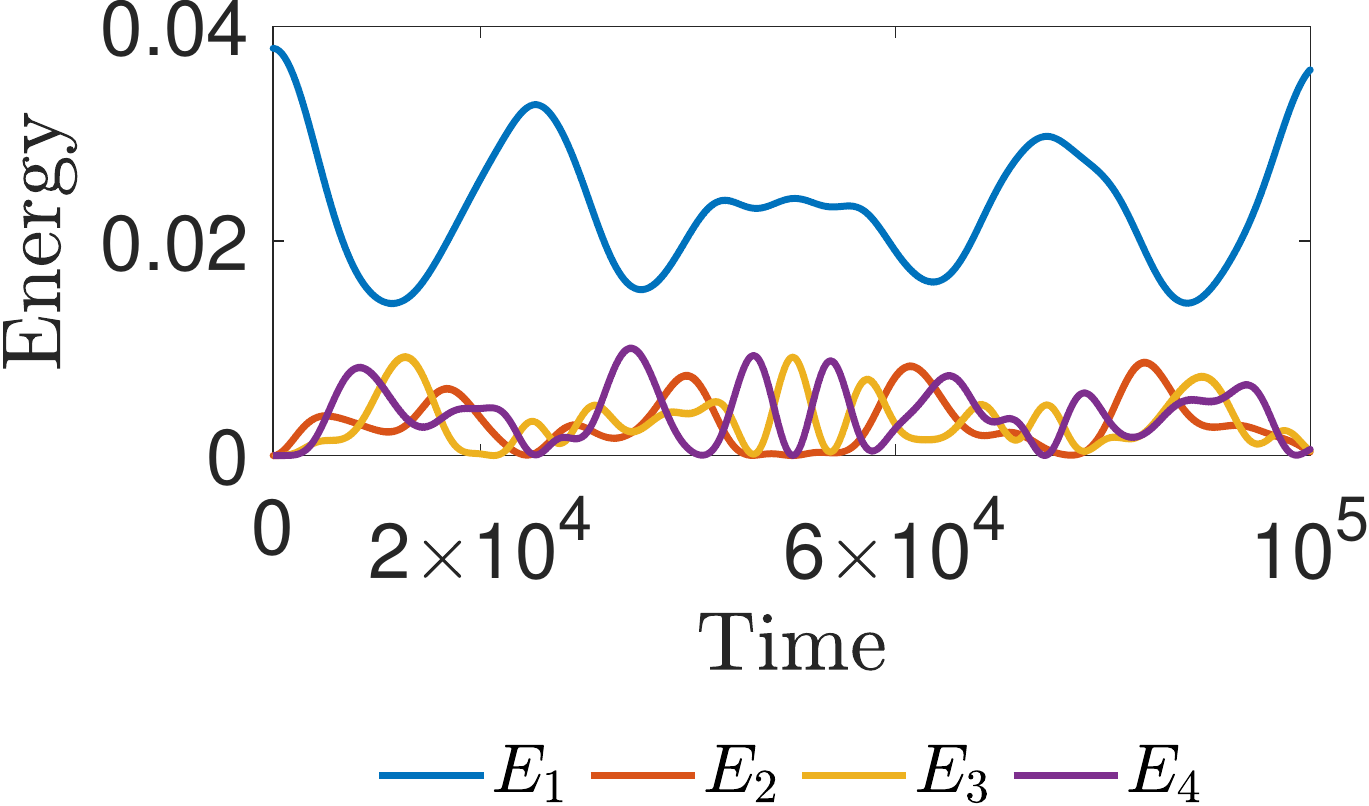}
		\caption{}
		\label{}
	\end{subfigure}%
	\begin{subfigure}{.3\textwidth}
		\centering
		\includegraphics[width=0.9\linewidth]{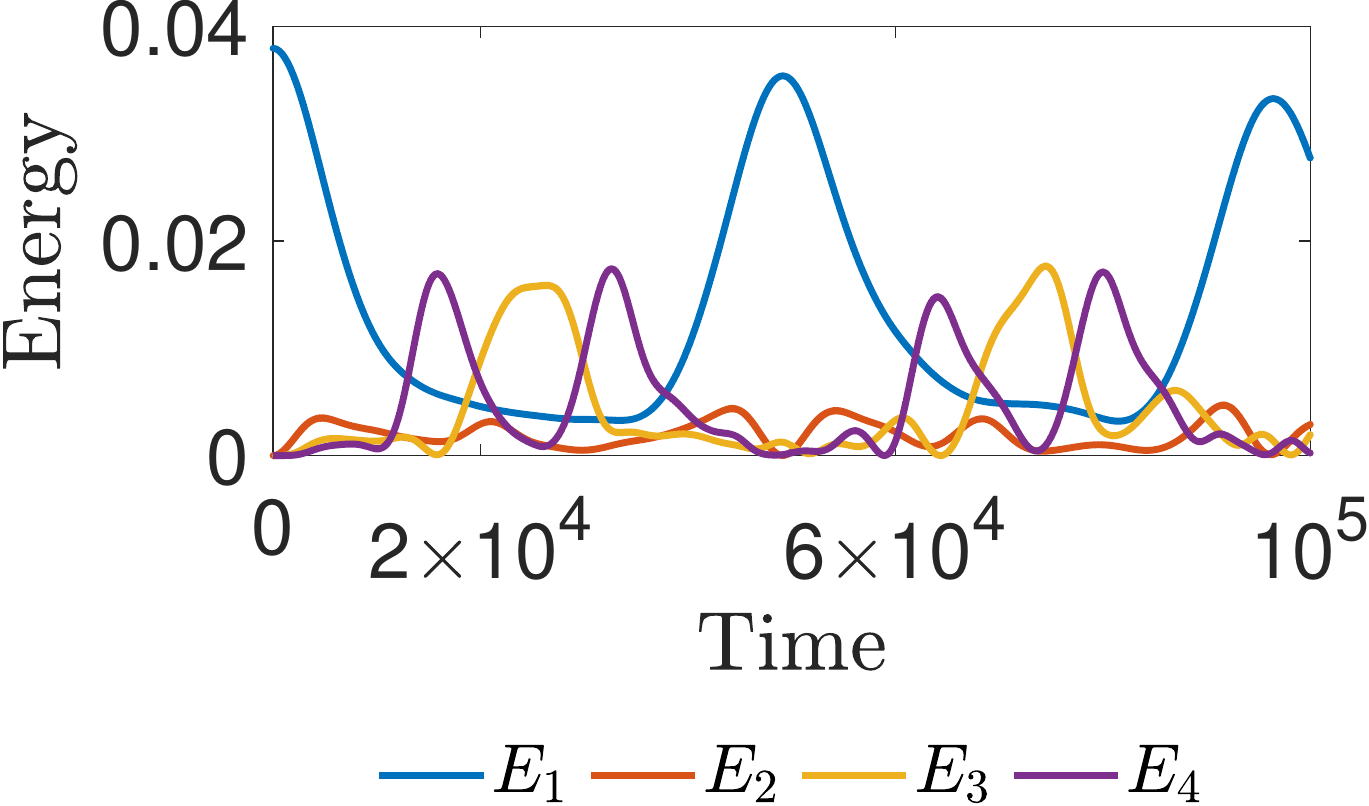}
		\caption{}
		\label{}
	\end{subfigure}\\
	\begin{subfigure}{.3\textwidth}
		\centering
		\includegraphics[width=0.9\linewidth]{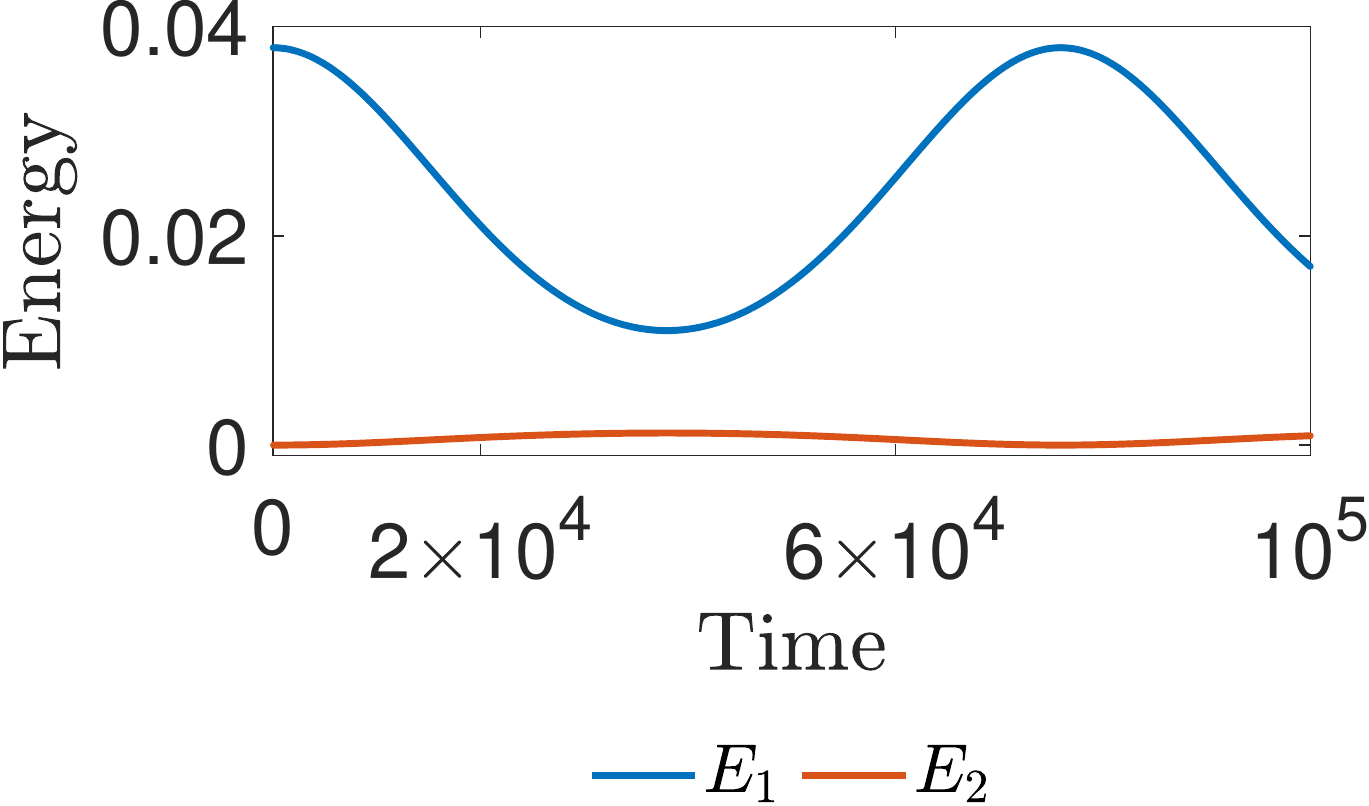}
		\caption{}
		\label{}
	\end{subfigure}%
	\begin{subfigure}{.3\textwidth}
		\centering
		\includegraphics[width=0.9\linewidth]{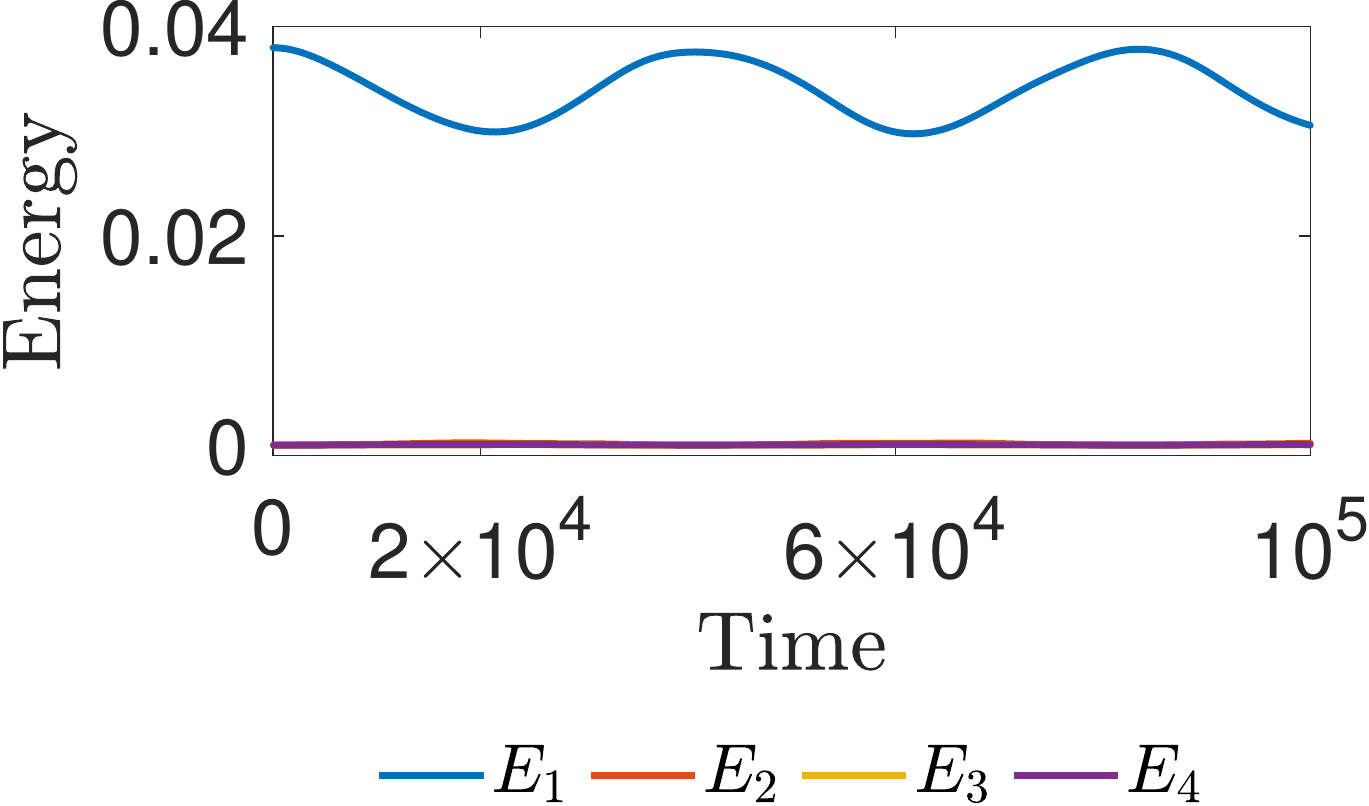}
		\caption{}
		\label{}
	\end{subfigure}%
	\begin{subfigure}{.3\textwidth}
		\centering
		\includegraphics[width=0.9\linewidth]{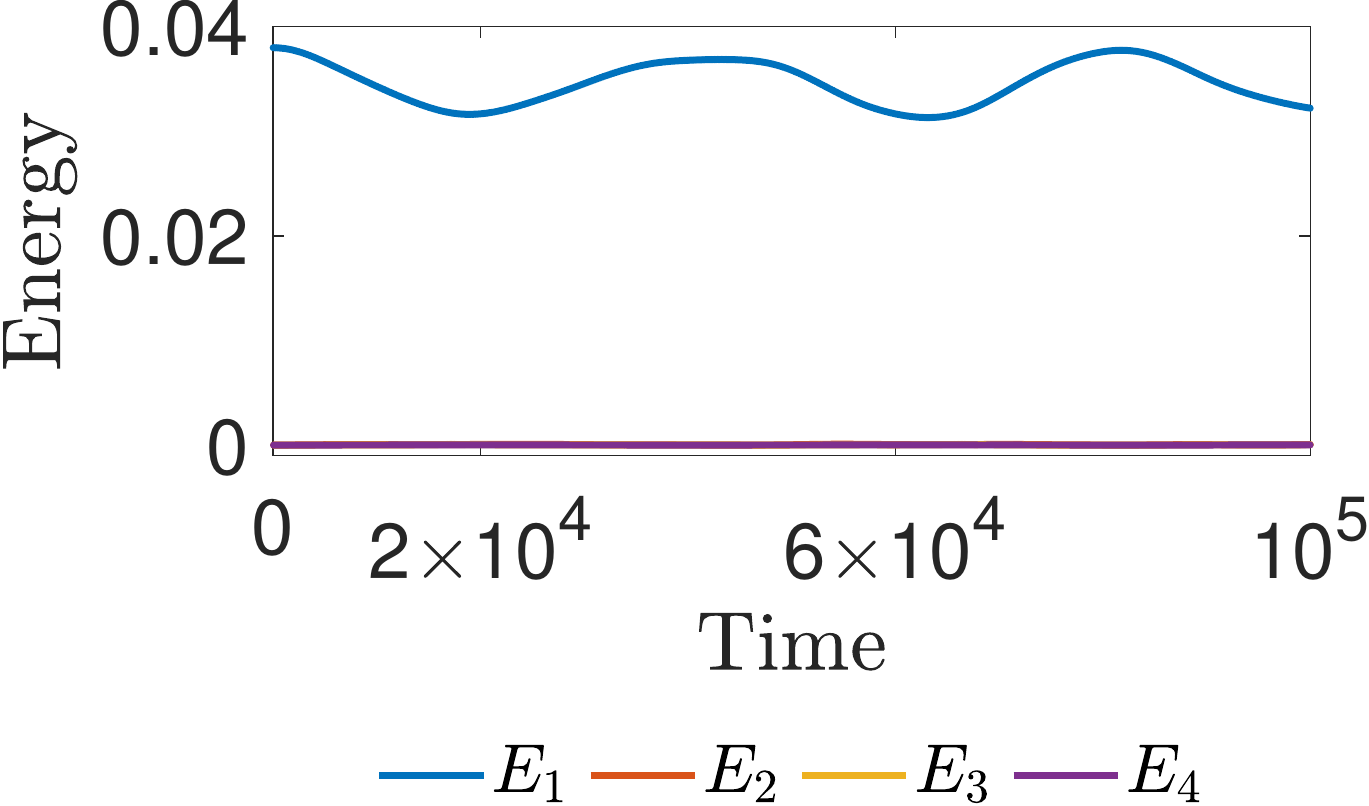}
		\caption{}
		\label{}
	\end{subfigure}
	\caption{Normal-mode energy in time obtained from integrating Eq.\ \eqref{eqnmfsvarnonlin} using only 2 normal modes in panels (a), (d), (g), 4 modes in panels (b), (e), (h) and 8 modes in panels (c), (f), (i). The tolerance is $0\%$ in panels (a) - (c), $5\%$ in panels (d) - (f), and $10\%$ in panels (g) - (i). We note that for illustration purposes, we plot in all panels only the normal mode energy of the first four modes and that all tolerances are smaller than $\tau_c$. Despite the fact that the last four modes are activated for $0\%$ and $5\%$ tolerance, they are essentially zero for $10\%$ tolerance in panel (i).}\label{fig:fig4}
\end{figure}

In the following, we show why the transfer of energy between modes reduces with the increase of the percentage of variability. This leads to a localized state in the energy-mode space and to finite-time blow up of solutions if variability is greater than $\tau_c$. When energy localization occurs, the plots of the normal mode energies suggest that most of the mode coordinates are vanishing in time. Therefore, we prefer to work in the normal-mode coordinate system than in the real (physical) space. In this framework, the equations of motion \eqref{varnonlin} can be written in the normal-mode coordinates space in a similar manner as in Eq.\ \eqref{eqnmotionsnmcoord}, namely in the form
\begin{align}
	\ddot{\mathbf{Q}} &=D\mathbf{Q}+A^{-1}\mathbf{\widehat{F}(Q)},
	\label{eqnmfsvarnonlin}
\end{align} 
for some nonlinear, vector-function $\mathbf{\widehat{F}(Q)}$ that depends on $\tau$, which is different to $\mathbf{F(Q)}$ in Eq.\ \eqref{eqnmotionsnmcoord} in the absence of variability. Our main assumption is that we can approximate system \eqref{eqnmfsvarnonlin} by considering only the first few modes. To illustrate numerically that this assumption is reasonable, we present in Fig.\ \ref{fig:fig4} the normal-mode energy for the set of equations of motion \eqref{eqnmfsvarnonlin} for 2, 4 and 8 normal modes and different percentages of variability. The parameter values in the set of equations of motion \eqref{eqnmfsvarnonlin} are calculated numerically for $N=64$ and the same percentage of variability as in Figs. \ref{fig:fig1} and \ref{fig:fig2}, where all remaining modes are set to 0 at all times.

Particularly, looking at Fig.\ \ref{fig:fig4}, we see that using $2$ and $4$ modes gives dynamics that are quantitatively different from those in Figs. \ref{fig:fig1} and \ref{fig:fig2}, with respect to the recurrence period. Nevertheless, even with only $2$ modes, we can still observe energy recurrence and localization for increasing percentage of variability. Therefore, in the following, we will consider a two normal-mode system in Eq.\ \eqref{eqnmfsvarnonlin}.

\section{A two normal-mode system and bifurcation analysis} \label{sec4}

Figure \ref{fig:fig4} suggests that when energy localization in the first few normal mode occurs, all higher modes have relatively much smaller energy. This gives us the idea that we can approximate Eq.\ \eqref{eqnmfsvarnonlin} by setting $Q_k(t)=0$ for $k=3,4,\dots,N$, and obtain the following two normal-mode system
\begin{subequations}
	\begin{align}
		\ddot{Q}_1 &= - \omega_1^2 Q_1 + \epsilon \left( A_1 Q_1^2+A_2Q_2^2+A_3Q_1Q_2 \right), \\
		\ddot{Q}_2 &= - \omega_2^2 Q_2 + \epsilon \left( B_1 Q_1^2+B_2Q_2^2+B_3Q_1Q_2 \right), 
	\end{align}
	\label{msvar}
\end{subequations}
where $A_i$, $B_i \in \mathbb{R}$, $i=1,2,3$ and $\omega_k$ is given in Eq.\ \eqref{wk}.

\subsection{Multiple-scale expansions}

Since $\omega_2 = 2 \omega_1 + \epsilon$, $|\epsilon| \ll 1$, we take the following asymptotic series 
\begin{subequations}
	\begin{align}
		Q_1 &=X_0(t,T) + \epsilon X_1(t,T) + \ldots,\\
		Q_2 &=Y_0(t,T) + \epsilon Y_1(t,T) + \ldots, 
	\end{align}
	\label{Q}
\end{subequations}
where $T=\epsilon t$ is a slow-time variable. The leading-order approximations to Eqs. \eqref{Q} are given by 
\begin{align}
	X_0 &= q_1(T) e^{i\omega_1 t} + q_1^*(T) e^{-i \omega_1 t}, \quad 
	Y_0 = q_2(T) e^{i\omega_2 t} + q_2^*(T) e^{-i \omega_2 t}. \label{msY03}
\end{align}

Substituting Eqs. \eqref{Q}, \eqref{msY03} into Eq.\ \eqref{msvar}, expanding the equations in $\epsilon$ and applying the standard solvability condition to avoid secular terms appearing (see e.g., \cite{syafwan2010discrete}), we obtain
\begin{subequations}
	\begin{align}
		i \frac{dq_1(T)}{d T} &= q_1(T) + \widetilde{A}q_1^*q_2,\\
		\quad i \frac{dq_2(T)}{d T} &= q_2(T) + \widetilde{B}q_1^2,
	\end{align}
	\label{msfinal2}
\end{subequations}

for $q_1$ and $q_2$, respectively, where $\tilde{A}=A_3/(2\omega_1)$ and $\tilde{B}=B_1/(2\omega_2)$. In this context, $i$ is the imaginary unit of the complex numbers. Following Eqs. \eqref{initxi}, the initial conditions of system \eqref{msfinal2} are given by
\begin{align}
	q_1(0)&=\frac{Q_1(0)}{2}, \label{ic-q1}\\ 
	q_2(0)&=0. \label{ic-q2}
\end{align}
We note that parameters $\tilde{A}$ and $\tilde{B}$ depend on $\tau$. 

\begin{figure}[h]
	\centering
	\begin{subfigure}{.5\textwidth}
		\centering
		\includegraphics[width=0.9\linewidth]{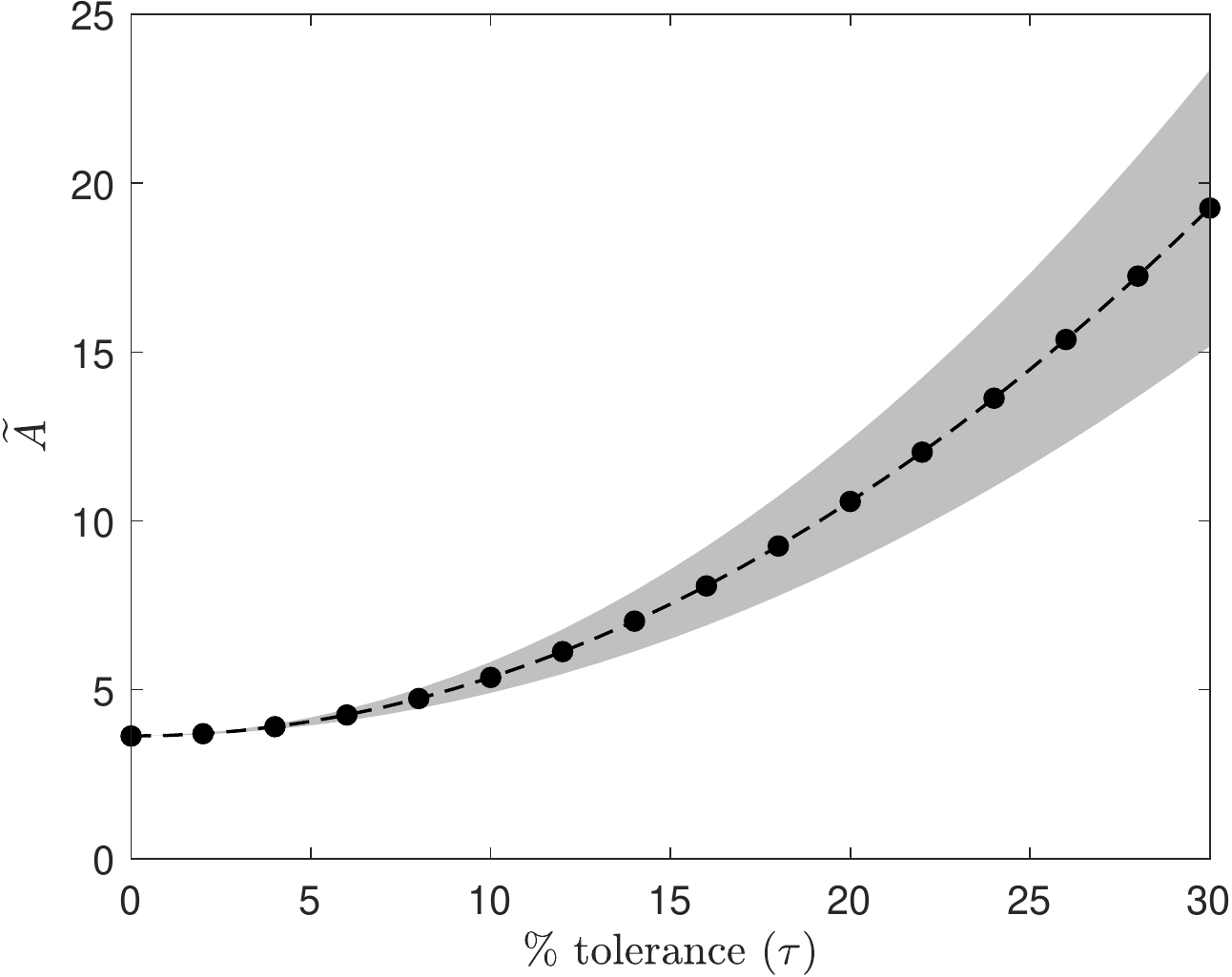}
		\caption{}
		\label{}
	\end{subfigure}%
	\begin{subfigure}{.5\textwidth}
		\centering
		\includegraphics[width=0.9\linewidth]{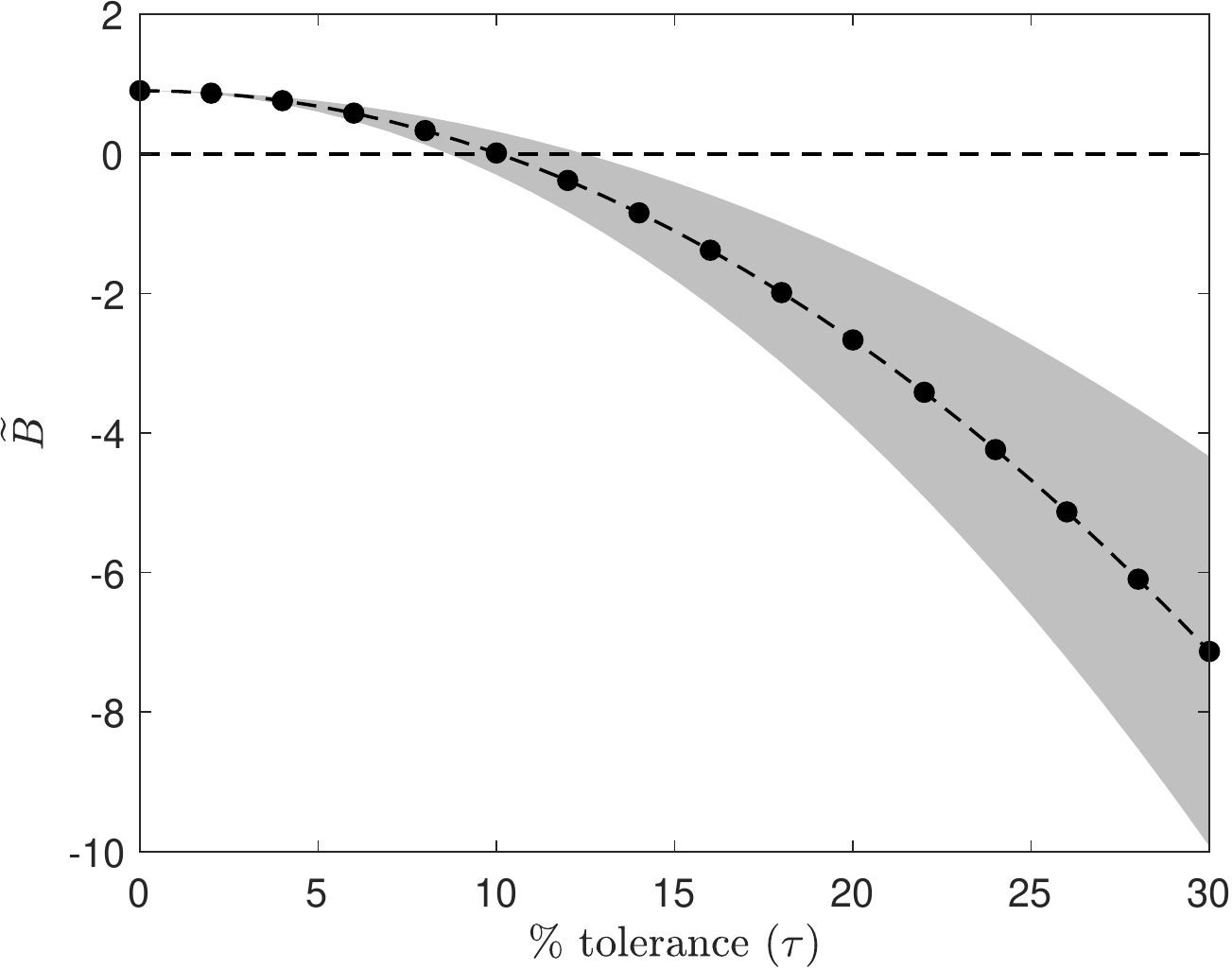}
		\caption{}
		\label{}
	\end{subfigure}%
	\caption{Plot of $\tilde{A}$ (in panel (a)) and $\tilde{B}$ (in panel (b)) as a function of the tolerance obtained numerically for $N=64$. The dash-dotted curve is the mean value over $100$ realisations of the same percentage of variability (see the discussion in the text), while the lengths of the shaded regions are two standard deviations. Using a polynomial regression, the mean is found to be given approximately by $\tilde{A}\approx0.01739\tau^2-0.00029\tau+3.62805$ and $\tilde{B}\approx-0.00893\tau^2-0.000084\tau+0.90728$, where the sums of square errors are $2.14\times10^{-15}$ and $3.46\times10^{-19}$ in panels (a) and (b), respectively.  Note the horizontal black dashed line at $\tilde{B}=0$ from which $\tau_c$ is derived (see text for details).}\label{fig:fig5}
\end{figure}

In Fig.\ \ref{fig:fig5}, we plot these parameters as a function of $\tau$ for $N=64$ particles and 100 realizations. These realizations have been computed by fixing $\tau$ and then opting for 100 sets of $N=64$ randomly generated numbers from the Gaussian distribution with mean 1 and standard deviation $\sigma=1/3 \times 0.01\tau$. Therefore, the $v_j$s in the 100 sets lie in the interval $[1-0.01\tau,1+0.01\tau]$. As we can see in panel (a), $\tilde{A}$ is positive for all $\tau$, whereas $\tilde{B}$ changes sign at around $\tau=10\%$. Particularly, $\tilde{B}$ starts positive for small $\tau$ values before it becomes negative at around $\tau=10\%$. By using polynomial regression, we have been able to fit the mean of the 100 realisations in panel (b) by the function $\tilde{B}\approx-0.00893\tau^2-0.000084\tau+0.90728$, with a sum of square errors (SSE) of $3.46\times10^{-19}$. This allowed us to estimate with good accuracy the threshold for the percentage of variability where $\tilde{B}$ changes sign and found to be given by $\tau_c\approx10.0749\%$ as $\tilde{B}(\tau_c)=0$. In Sec.\ \ref{equilibrium_solutions}, we show that when $\tilde{B}<0$, that is for $\tau>\tau_c$, trajectories of Eqs. \eqref{eqnmotions} may blow up in finite time.

\begin{figure}[h]
	\centering
	\begin{subfigure}{.5\textwidth}
		\centering
		\includegraphics[width=0.9\linewidth]{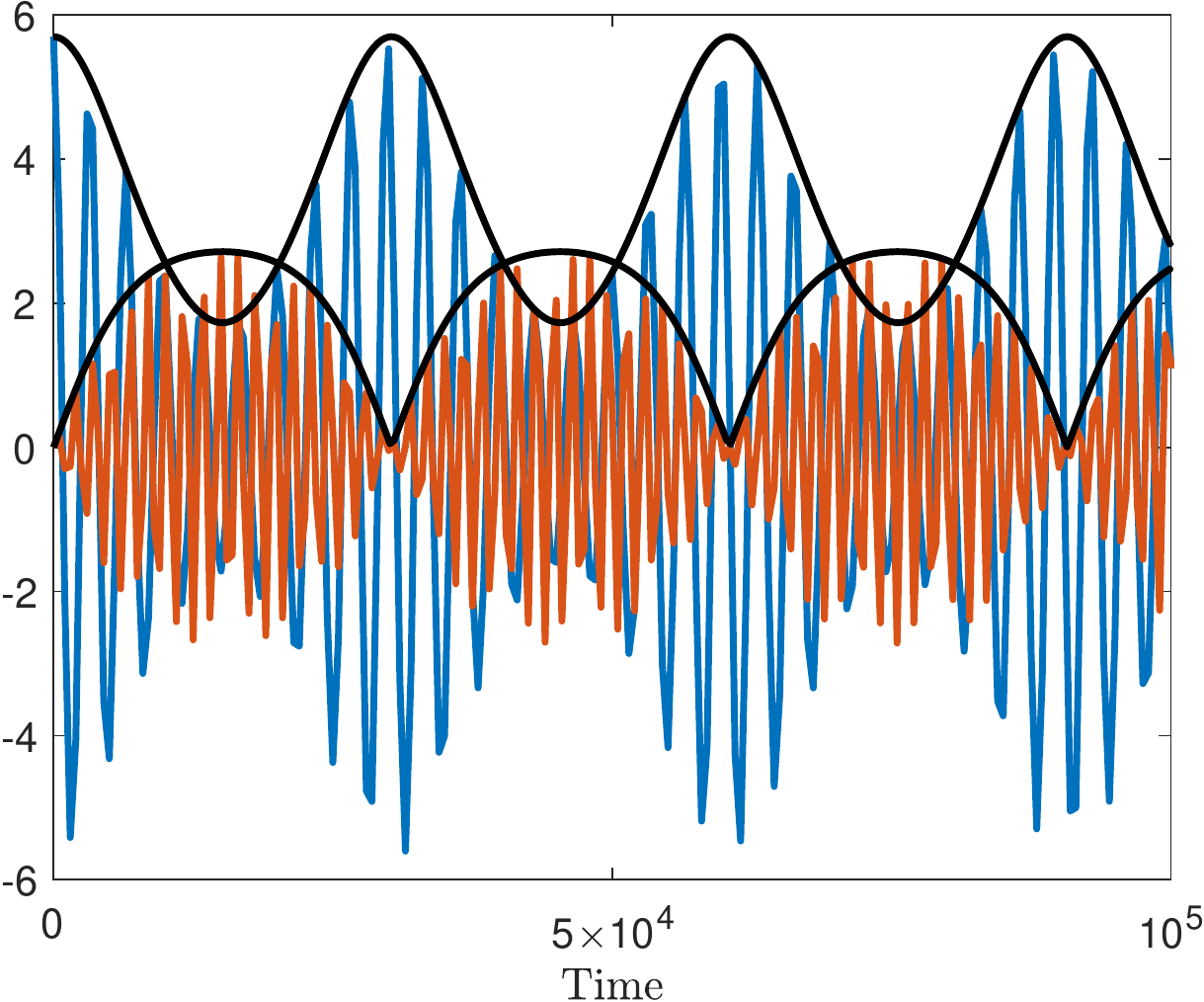}
		\caption{}
		\label{}
	\end{subfigure}%
	\begin{subfigure}{.5\textwidth}
		\centering
		\includegraphics[width=0.9\linewidth]{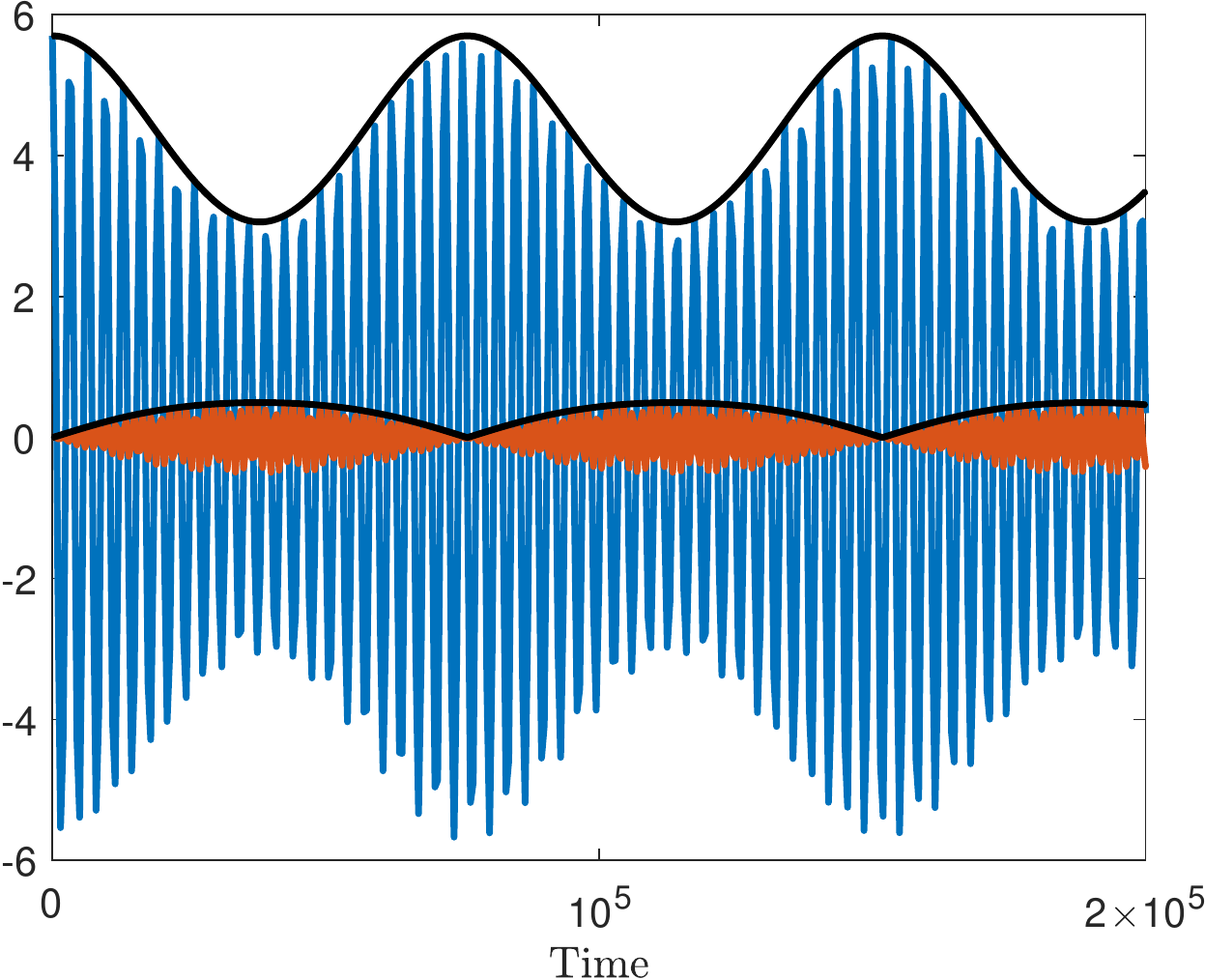}
		\caption{}
		\label{}
	\end{subfigure}
	\caption{Time evolution of the normal mode variables $Q_1$ (blue curve) and $Q_2$ (red curve) with their envelopes $q_1$ and $q_2$ (black curves) from Eqs. \eqref{msfinal2} for $\tau=0\%$ in panel (a) and $\tau=10\%$ in panel (b). Note that in both panels $\tau<\tau_c$, so trajectories do not blow up.}\label{fig:fig6}
\end{figure}

A comparison of the dynamics of the normal modes $Q_1$ and $Q_2$ of Eq.\ \eqref{msvar} and those of the slow-time variables $q_1$ and $q_2$ of Eqs. \eqref{msfinal2} is shown in Fig.\ \ref{fig:fig6}, where one can see that $q_j$ is an envelope of $Q_j$ for $j=1,2$.

Next we explain the cause of localization with the increase of the percentage of variability $\tau$. Note that from Eqs. \eqref{msfinal2}, there can be transfer of energy from $q_1(t)$ to $q_2(t)$ through the nonlinear coupling coefficient $\tilde{B}$. Panel (b) in Fig.\ \ref{fig:fig5} shows that $\tilde{B}$ decreases from positive values with the increase of $\tau$ until $\tau=\tau_c$, after which it becomes negative. When $\tilde{B}$ vanishes at $\tau=\tau_c$, there is no transfer of energy and hence localization. In the following, we will also show that when $\tilde{B}<0$, i.e., for $\tau>\tau_c$, there might be unbounded trajectories that blow up in finite time.

\subsection{Equilibrium solutions}\label{equilibrium_solutions}

We start by analyzing the standing wave solutions of the envelope equations \eqref{msfinal2}. To do so, it is convenient to write $q_1$ and $q_2$ in polar form $q_1=r_1 e^{i\phi_1}$ and $q_2=r_2 e^{2i\phi_2}$, where $r_1=|q_1|$, $r_2=|q_2|$. Then, we define the new variables 
\begin{subequations}
	\begin{align} 
		P&=r_1^2+r_2^2, \label{P} \\
		\Delta &= r_1^2-r_2^2, \label{delta}\\
		\theta &= \phi_2 - \phi_1. \label{theta}
	\end{align}
	\label{transf}
\end{subequations}
These variables satisfy the set of equations (see \cite{pickton2013integrability} for a similar derivation)
\begin{subequations}
	\begin{align}
		\dot{P} &= \frac{\widetilde{A}-\widetilde{B}}{\widetilde{A}+\widetilde{B}} \dot{\Delta}, \label{P2}\\
		\dot{\Delta}&=\frac{\sqrt {2\left(P-\Delta\right)}\sin
			\left( 2\,\theta \right)  \left( P+\Delta \right)  \left( \widetilde{A}+\widetilde{B}
			\right)}{2}, \label{delta3} \\
		\dot{\theta}&=-{\frac {2\,\widetilde{A}\cos \left( 2\,\theta
				\right) (\Delta-P)+\widetilde{B}\cos \left( 2\,
				\theta \right) (\Delta+P)-\sqrt {2\left(P-\Delta\right)}}{2\sqrt {2\left(P-\Delta\right)}}}.\label{theta3}
	\end{align}
	\label{reduced}
\end{subequations}

\begin{figure}[h]
	\centering
	\includegraphics[width=0.5\textwidth]{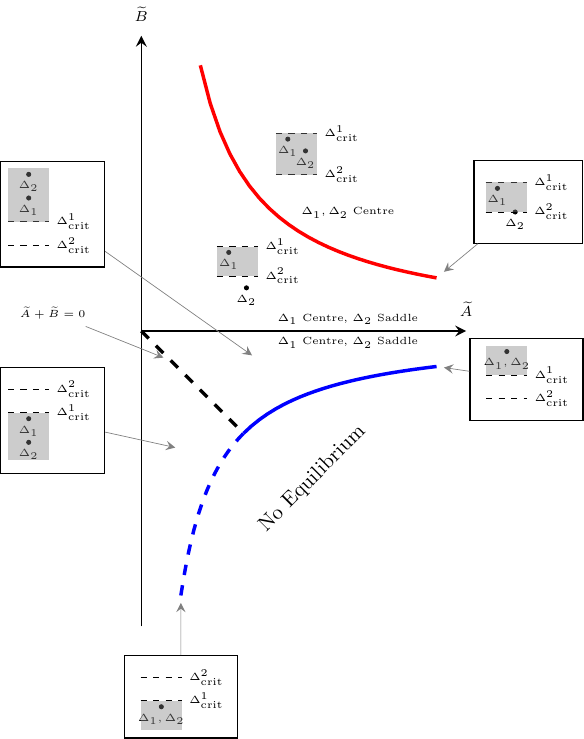}%
	\caption{Bifurcation diagram of the equilibrium points $\Delta_1$ and $\Delta_2$ and the regions where the dynamics of system \eqref{reduced} is well-defined (see text for more details).}
	\label{fig:fig7}%
\end{figure}

From Eq.\ \eqref{P2}, the constant of motion $C$ follows
\begin{align}
	C&=P-\frac{\widetilde{A}-\widetilde{B}}{\widetilde{A}+\widetilde{B}} \Delta. \nonumber
\end{align}
System \eqref{msfinal2} is transformed into Eqs. \eqref{reduced} by using Eqs. \eqref{transf}, where we assume that $r_1$ and $r_2$ are non negative real numbers. Equation \eqref{reduced} requires $P-\Delta>0$ in order to have real-valued solutions, whereas Eq.\ \eqref{transf} requires $P-\Delta>0$ and $P+\Delta\geq0$, otherwise $r_1$ and $r_2$ will be complex numbers. We call the region which satisfies these two inequalities the well-defined region and denote it by the shaded area in Fig.\ \ref{fig:fig7}. $\Delta_2$ is outside the shaded region in the area below the red curve and above $\widetilde{B}=0$. This implies that $\Delta_2$ is the equilibrium of system \eqref{reduced} only, but not of system \eqref{msfinal2}.

The latter result implies that the dynamics of Eq.\ \eqref{msfinal2} can be described by the remaining equations \eqref{delta3} and \eqref{theta3}, i.e., in terms of $\Delta$ and $\theta$ only. As discussed before, Eqs. \eqref{delta3}, \eqref{theta3} are valid only when $P-\Delta > 0$. However, Eqs. \eqref{P}, \eqref{delta} imply that $P+\Delta \geq 0$. These two inequalities determine the region where Eqs. \eqref{P}, \eqref{delta} are defined in the $(\Delta,\theta)$-plane. As this region depends on $\widetilde{A}$ and $\widetilde{B}$, we consider the following cases:

\begin{itemize}
	\item{If $\frac{\widetilde{A}-\widetilde{B}}{\widetilde{A}+\widetilde{B}} \geq 1$, then $\Delta > \text{max}\left\{\Delta_{\text{crit}}^{1},\Delta_{\text{crit}}^{2}\right\}$,
		where
		\begin{align*}
			\Delta_{\text{crit}}^{1}=\frac{C\left(\widetilde{A}+\widetilde{B}\right)}{2\widetilde{B}}\mbox{ and }\Delta_{\text{crit}}^{2}=-\frac{C\left(\widetilde{A}+\widetilde{B}\right)}{2\widetilde{A}}.
		\end{align*} 
	}
	\item{If $-1 \leq \frac{\widetilde{A}-\widetilde{B}}{\widetilde{A}+\widetilde{B}} < 1$, then $\Delta_{\text{crit}}^{2}\leq\Delta < \Delta_{\text{crit}}^{1}$. 
	}
	\item{If $\frac{\widetilde{A}-\widetilde{B}}{\widetilde{A}+\widetilde{B}} < -1$, then $\Delta < \text{min}\left\{\Delta_{\text{crit}}^{1},\Delta_{\text{crit}}^{2}\right\}$.}
\end{itemize}

The regions, in which the reduced system \eqref{reduced} is well-defined, are plotted in the ($\tilde{A}$,$\tilde{B}$)-plane in Fig.\ \ref{fig:fig7}. 

To study the reduced system of Eqs. \eqref{delta3}, \eqref{theta3}, we restrict the phase difference $\theta$ in the interval $0\leq \theta < \pi$ and obtain two equilibrium points, namely $(\theta_j,\Delta_j)$, $j=1,2$, where
\begin{equation}\label{theta_1_2}
	\theta_1=0\mbox{ or } \pi/2\mbox{ and }\theta_2= \pi/2.
\end{equation}
Particularly, there are two cases with respect to $\tilde{B}$. The first one is when $\widetilde{B}>0$, in which case $\theta_1=0$ or $\theta_2=\pi/2$, and the second when $\widetilde{B}<0$, in which case $\theta_1=\theta_2=\pi/2$. Then, 
\begin{align}
	\Delta_{j}&={\frac { \left( 6\,{\widetilde{A}}^{2}C-3\,\widetilde{A}\widetilde{B}C-(-1)^j\sqrt {1+6\,\widetilde{A}
				\left( \widetilde{A}+\widetilde{B} \right) C}-1 \right)  \left( \widetilde{A}+\widetilde{B} \right) }{18{\widetilde{A}}^{2}\widetilde{B}}} \label{deltaequi}.
\end{align}
The stability of the equilibrium points is determined by the eigenvalues of the Jacobian matrix of Eqs. \eqref{delta3}, \eqref{theta3}, evaluated at the equilibrium points, i.e., by
\begin{align}
	\lambda^{(j)}_{1,2} &= \pm\frac{\sqrt {-3-18{\widetilde{A}}^{2}C-18\widetilde{A}\widetilde{B}C+6(-1)^j\sqrt {1+6\,\widetilde{A}\left( \widetilde{A}+\widetilde{B} \right) C}}}{3}. \label{eigequi}
\end{align}

From Eq. (\ref{deltaequi}) it follows that the equilibrium points exist when 
\begin{align}
	1+6\,\widetilde{A}
	\left( \widetilde{A}+\widetilde{B} \right) C \geq 0. \label{equicrit}
\end{align}
For the initial conditions (\ref{ic-q1}), (\ref{ic-q2}), Eq. (\ref{equicrit}) becomes
\begin{align}
	1+12\widetilde{A}\widetilde{B}r_1^2 \geq 0. \nonumber
\end{align}
Therefore, the threshold for the existence of the equilibrium is given by
\begin{align}
	1+12\widetilde{A}\widetilde{B}r_1^2= 0,\nonumber
\end{align}
which is the blue curve in Fig.\ \ref{fig:fig7}. The dashed and solid lines represent the curve below and above the line $\widetilde{A}+\widetilde{B}=0$, respectively.

Comparing Eqs. (\ref{deltaequi}) and (\ref{eigequi}), we conclude that when the equilibrium points exist, they are either a centre or a saddle node. Particularly, for the initial conditions in Eqs. \eqref{ic-q1}, \eqref{ic-q2}, the thresholds for the eigenvalues that discriminate between a centre and a saddle node are
\begin{align}
	1+12\widetilde{A}\widetilde{B}r_1^2 &= 0, \label{eigcrit1} \\
	1-4\widetilde{A}\widetilde{B}r_1^2 &= 0.  \label{eigcrit2}
\end{align} 
In Fig. (\ref{fig:fig7}), we plot $\tilde{A}+\tilde{B}=0$, Eqs. \eqref{eigcrit1} and \eqref{eigcrit2} as the black dashed, blue and red curves, respectively.

System (\ref{reduced}) with parameter values above the red curve in Fig. (\ref{fig:fig7}) is bounded, with $\Delta_{\text{crit}}^1$ and $\Delta_{\text{crit}}^2$ being the upper and lower bounds, respectively. The two equilibria given in Eq.\ \eqref{deltaequi} are both centres, and are therefore stable. When the parameter values lie on the red curve, $\Delta_2=\Delta_{\text{crit}}^2$. Furthermore, if the parameter values are below the red curve and $\widetilde{B}>0$, system \eqref{reduced} is still bounded, but it only shares one equilibrium point $\Delta_1$ with system \eqref{msfinal2}, whereas $\Delta_2$ does not belong to the well-defined region. Equation \eqref{reduced}, on the other hand, is unbounded when $\widetilde{B}<0$. In this case, it either extends to $\Delta\to\infty$ or $-\infty$ and depending on the value of $\frac{\widetilde{A}-\widetilde{B}}{\widetilde{A}+\widetilde{B}}$, $\Delta_1$ can be a centre and $\Delta_2$ a saddle node in this region. Additionally, the system has only one equilibrium on the blue curve.

\begin{figure}[h]
	\centering
	\begin{subfigure}{.5\textwidth}
		\centering
		\includegraphics[width=0.9\linewidth]{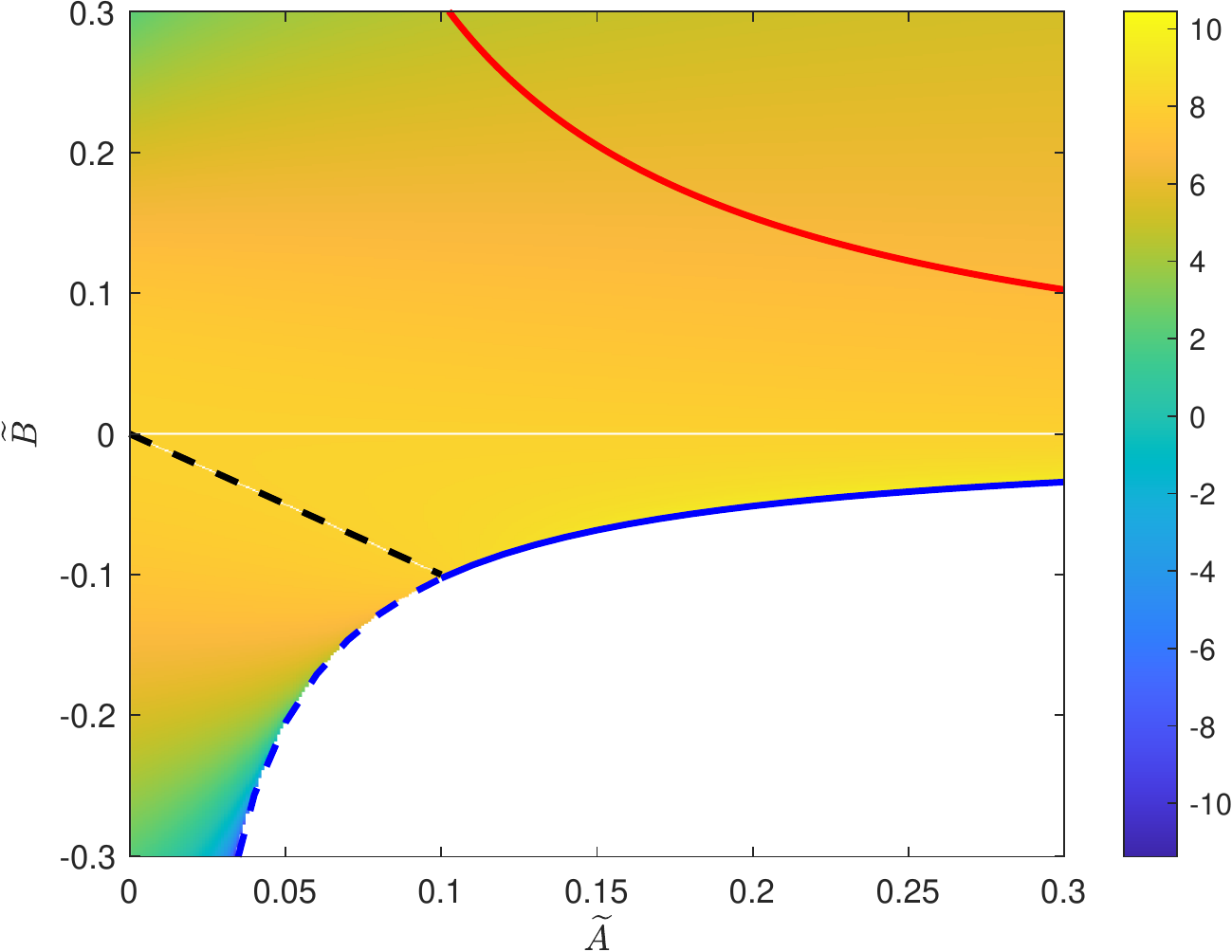}
		\caption{}
		\label{fig8a}
	\end{subfigure}%
	\begin{subfigure}{.5\textwidth}
		\centering
		\includegraphics[width=0.9\linewidth]{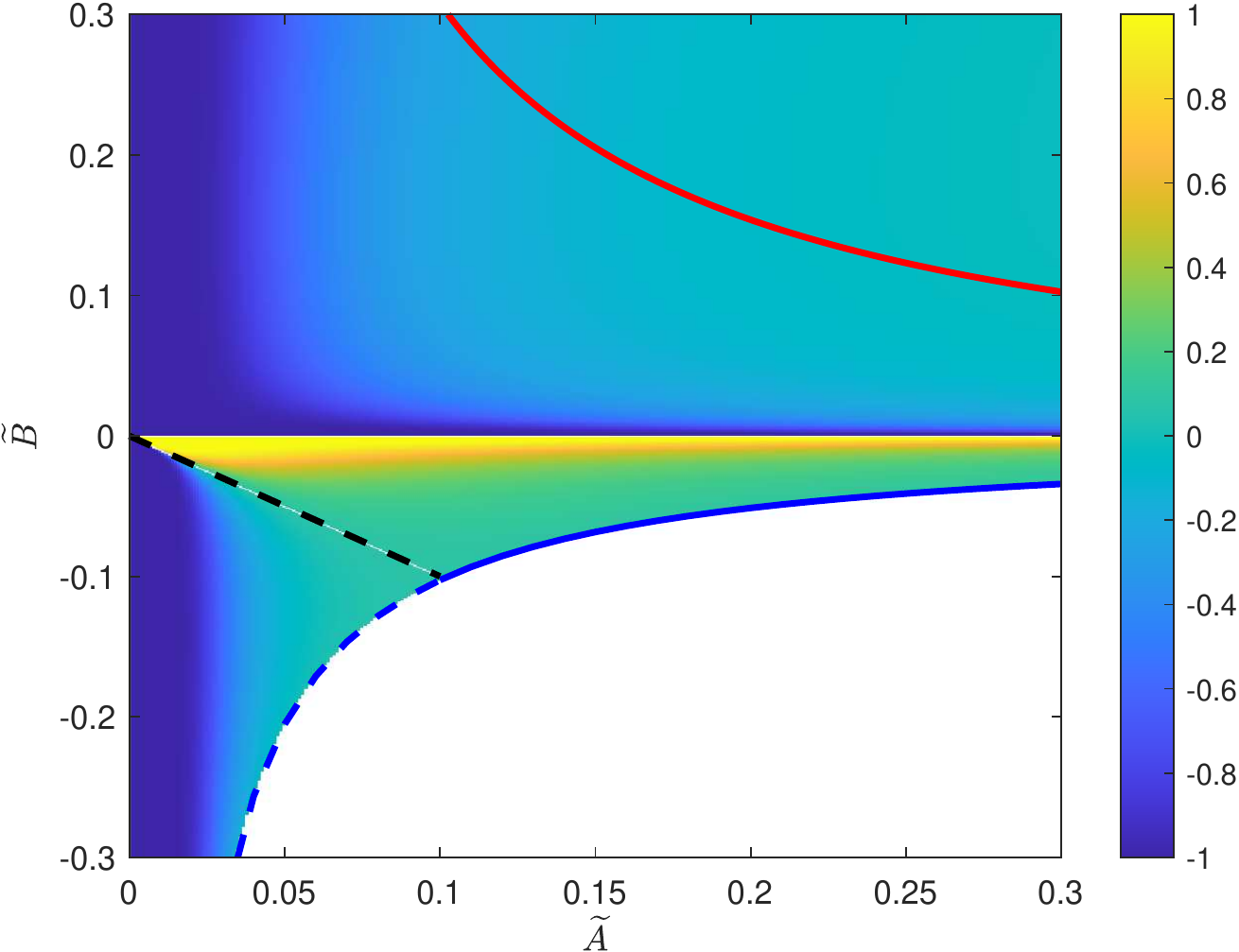}
		\caption{}
		\label{fig8b}
	\end{subfigure}
	\caption{Plot of (a) $\Delta_1$ and (b) $\tanh(\Delta_2/100)$ as a function of $\widetilde{A}$ and $\widetilde{B}$. In panel (a), the color bar denotes the values of $\Delta_1$ and in panel (b), the values of $\tanh(\Delta_{2}/100)$. The black dashed, red and blue curves are discussed in the text and are the same with those in Fig.\ \ref{fig:fig7}.}
	\label{fig:fig8}
\end{figure}

The location of the equilibrium points $(\theta_j,\Delta_j)$ in Eqs. \eqref{theta_1_2}, \eqref{deltaequi} and their nature are shown in Fig.\ \ref{fig:fig7}. We also plot the values of $\Delta_{j}$ in Fig.\ \ref{fig:fig8}. To better visualise $\Delta_{2}$ as it approaches infinity when $\widetilde{A}$ or $\widetilde{B}$ approaches zero, we plot in Fig.\ \ref{fig:fig8} (b) $\tanh(\Delta_{2}/100)$ instead of $\Delta_{2}$.

\begin{figure}[h]
	\centering
	\begin{subfigure}{.5\textwidth}
		\centering
		\includegraphics[width=\linewidth]{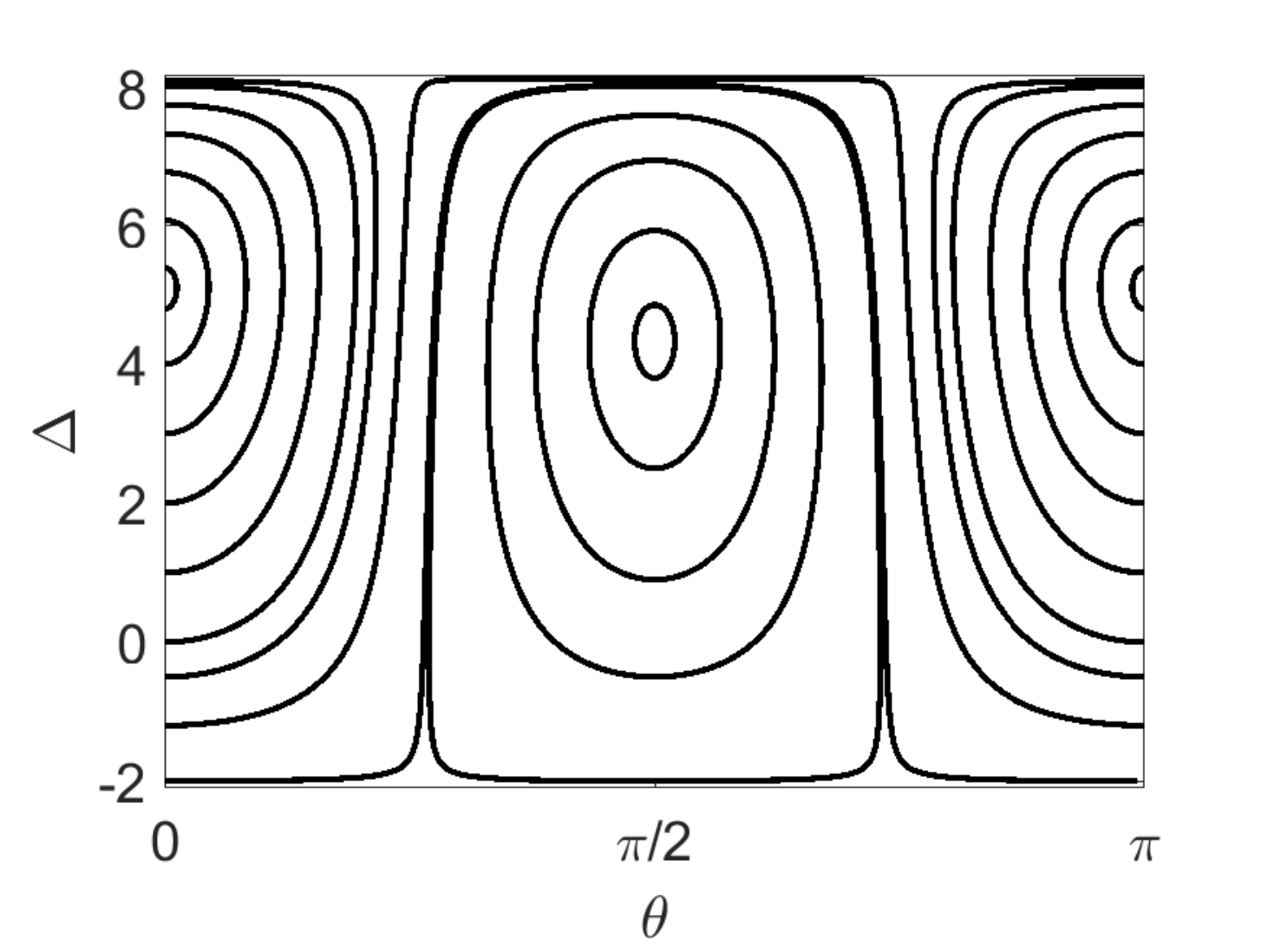}
		\caption{}
		\label{}
	\end{subfigure}%
	\begin{subfigure}{.5\textwidth}
		\centering
		\includegraphics[width=\linewidth]{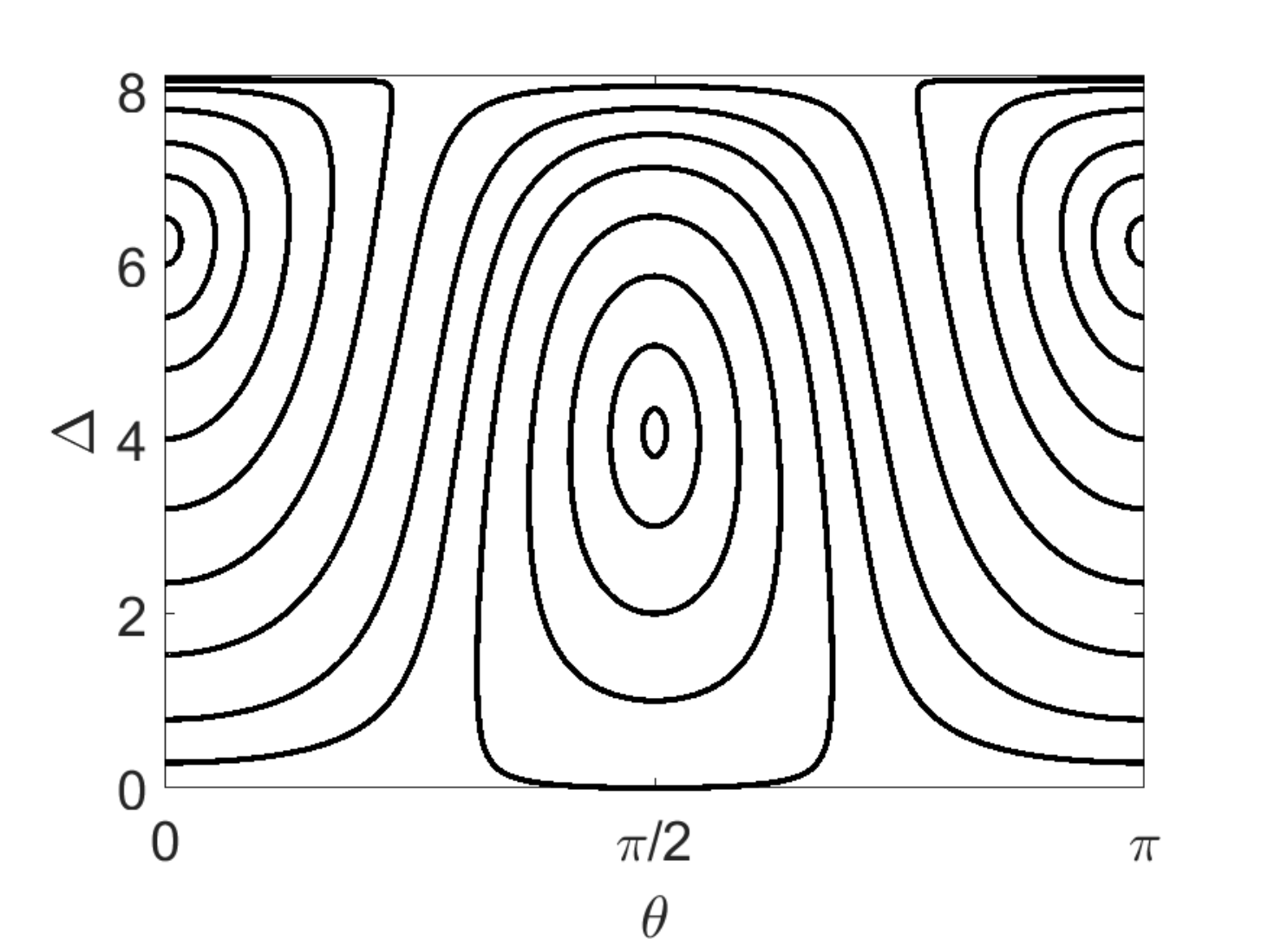}
		\caption{}
		\label{}
	\end{subfigure}
	\caption{Phase portraits of the reduced system of Eqs. \eqref{delta3}, \eqref{theta3} for (a) $\tau=0\%$ percentage of variability and (b) $\tau=10\%$ percentage of variability. }
	\label{fig:fig9}
\end{figure}

In the following, we illustrate the phase portrait of the reduced system of Eqs. \eqref{delta3}, \eqref{theta3} for different percentages of variability $\tau$, which correspond to different values of $\tilde{A}$ and $\tilde{B}$. When there is no variability (i.e., for $\tau=0\%$), the parameter values are $\widetilde{A}=3.63$ and $\widetilde{B}=0.91$ and the equilibrium points are $(\theta_1,\Delta_1)=(0,5.09)$ and $(\theta_2,\Delta_2)=(\pi/2,4.34)$. Both are stable and the phase space in this case is shown in Fig.\ \ref{fig:fig9}(a). As we can see in panel (b) in Fig.\ \ref{fig:fig5}, as $\tau$ increases, $\widetilde{B}$ decreases and becomes negative for $\tau>\tau_c$. The parameter values for $\tau=10\%$ variability are $\widetilde{A}=4.97$ and $\widetilde{B}=0.05$ and the equilibrium points are $(\theta_1,\Delta_1)=(0,6.27)$ and $(\theta_2,\Delta_2)=(\pi/2,4.08)$. Similar to the previous case, both equilibrium points are stable and the phase space is shown in Fig.\ \ref{fig:fig9}(b). Note that for the initial conditions \eqref{ic-q2}, we have that
\begin{align*}
	\lim_{\widetilde{B} \rightarrow 0} \Delta_{\text{crit}}^{1} &= r_1^2,\quad 
	\lim_{\widetilde{B} \rightarrow 0} \Delta_{\text{crit}}^{2} = 0,
\end{align*}
which shows that $\Delta$ becomes positive as we increase $\tau$. Indeed, $\Delta>0$ corresponds to energy localization as the magnitude of $q_1$ remains larger than $q_2$.

\begin{figure}[h]
	\centering
	\includegraphics[width=0.5\textwidth]{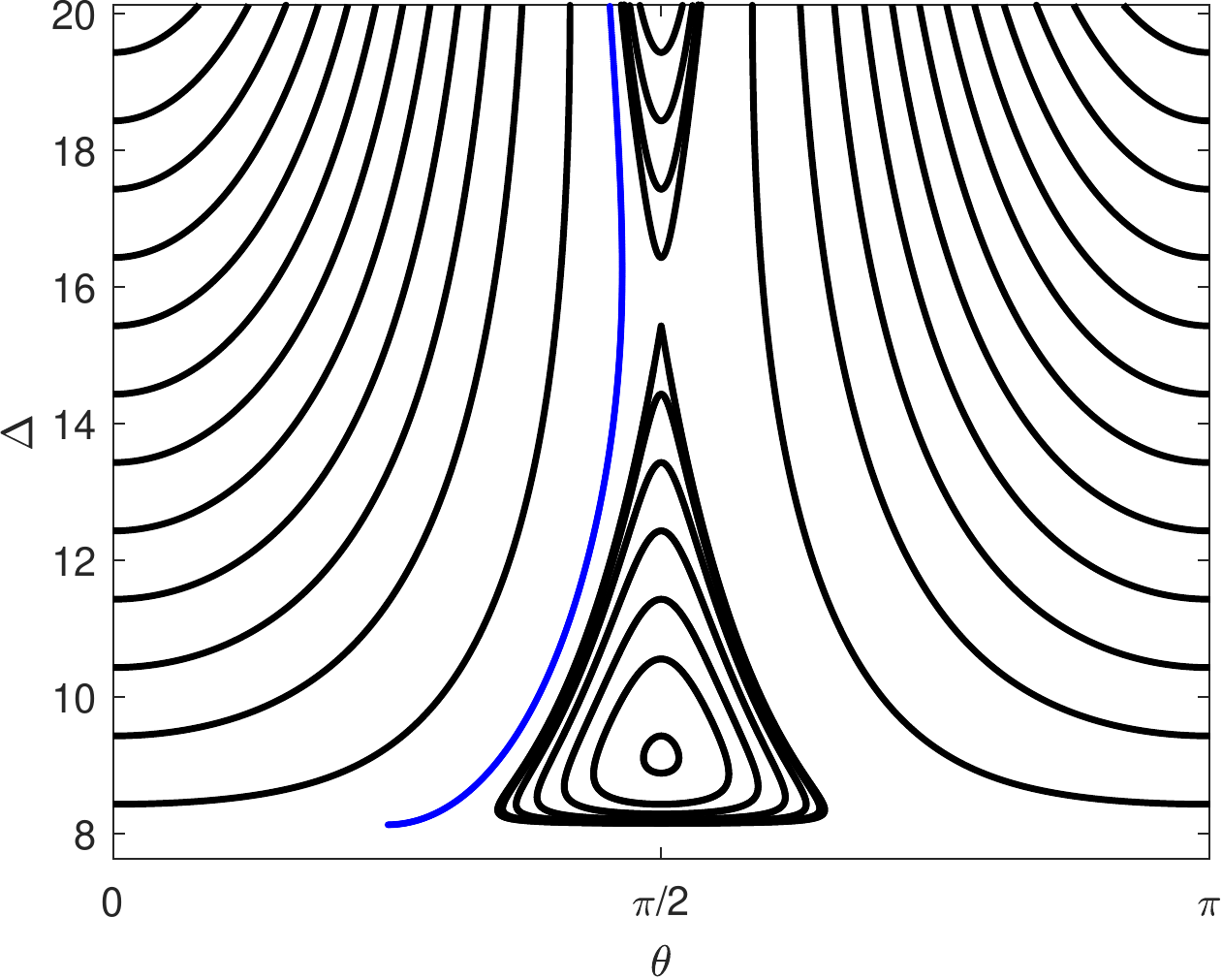}
	\caption{The same as Fig.\ \ref{fig:fig9}, where the parameter values are $\widetilde{A}=5.3932$ and $\widetilde{B}=-0.0015$, which correspond to $\tau\approx 10.0833\%>\tau_c$. The blue curve is the trajectory of the initial condition in Eqs. \eqref{ic-q1}, \eqref{ic-q2}.}
	\label{fig:fig10}
\end{figure}

As we can see in Fig.\ \ref{fig:fig10} for $\tau\approx 10.0833\%>\tau_c$, $\widetilde{B}$ is negative ($\widetilde{B}=-0.0015$) and the region in the $(\Delta,\theta)$-space becomes unbounded (see also Fig.\ \ref{fig:fig7}). It extends to either $\Delta\to\infty$ or $-\infty$ and depends on $\widetilde{A}$. In this case, the two equilibrium points are $(\theta_1,\Delta_1)=(\pi/2,9.1274)$, which is a (stable) center, and $(\theta_2,\Delta_2)=(\pi/2,15.4383)$, which is a (unstable) saddle point. The plot shows that in this case, one may obtain bounded solutions as well as unbounded ones, depending on the initial condition. For example, the initial condition of Eqs. \eqref{initxi} (or Eqs. \eqref{ic-q1}, \eqref{ic-q2}) results in $\theta$ and $\Delta$ values in the unbounded region in Fig.\ \ref{fig:fig10}, where the trajectory is shown as the blue curve and starts at the bottom of the plot.

\section{Chaotic behavior}\label{sec5}

Energy recurrences arise in the homogeneous FPUT lattice \eqref{Hamilt} when the system remains in the quasi-stationary state for an extremely long time, making the approach to equipartition of energy unobservable. In the quasi-stationary state, the FPUT lattice can be viewed as the perturbation of the regular, integrable Toda lattice \cite{ponno2011two}.

Here we study the effect of variability on the chaotic properties of system \eqref{varnonlin}. Particularly, we consider lattices of $N=4,8,16,32,64$ particles in systems \eqref{Hamilt} (homogeneous, no variability) and \eqref{varnonlin} (with variability) and use the maximum Lyapunov exponent (mLE) \cite{benettin1980lyapunov} and Smaller Alignment Index (SALI) \cite{skokos2003does,skokos2004detecting} to discriminate between regular and chaotic dynamics. We want to see if energy localization in the first normal mode that we observed in Secs. \ref{sec3} and \ref{sec4} for $\tau=10\%<\tau_c$ corresponds to chaotic dynamics, by increasing $\tau$ from 0 to $10\%$.

To compute mLE, we follow the evolution of a trajectory starting at the initial point $$\mathbf{x}(0)=(q_1(0),\ldots,q_N(0),p_1(0),\ldots,p_N(0)),$$ that evolves according to Hamilton's equations of motion
\begin{equation}
	\dot{\mathbf{x}}=\mathbf{f(x)}=\left[ \frac{\partial{H}}{\partial{\mathbf{p}}} \quad -\frac{\partial{H}}{\partial{\mathbf{q}}} \right]^T,
	\nonumber
\end{equation}
and the evolution of a deviation vector $$\mathbf{w}(0)=(\delta q_1(0),\ldots,\delta q_N(0),\delta p_1(0),\ldots,\delta p_N(0)),$$ that evolves according to the variational equation
\begin{align}
	\dot{\mathbf{w}}& =  \frac{\partial{\mathbf{f}}}{\partial{\mathbf{x}}}(\mathbf{x}(t)) \cdot \mathbf{w}.
	\label{wode}
\end{align}
Then mLE is defined as
\begin{align}
	\lambda= \lim_{t \to \infty}\frac{1}{t}\ln \frac{||\mathbf{w}(t)||}{||\mathbf{w}(0)||}, \nonumber
\end{align}
where $\ln$ is the natural logarithm. If mLE converges to zero following the law $1/t$, then the trajectory is regular, whereas if it converges to a positive value in time, then the trajectory is chaotic \cite{skokos2010lyapunov}. Hence it is convenient to plot mLE in $\log_{10}$$-\log_{10}$ scales as the law $1/t$ becomes then a line with negative slope and serves as a guide to the eye.

To compute SALI, we follow the evolution of the same initial condition and two deviation vectors $\mathbf{w}_1(0)$, $\mathbf{w}_2(0)$. Then, SALI is defined by
\begin{align*}
	\textrm{SALI}(t)=\min \{ \left\| \hat{\mathbf{w}}_1(t) - \hat{\mathbf{w}}_2(t)\right\|,\left\| \hat{\mathbf{w}}_1(t) + \hat{\mathbf{w}}_2(t)\right\| \},	
\end{align*}
where $\hat{\mathbf{w}}_i(t)= \frac{\mathbf{w}_i(t)}{\left\| \mathbf{w}_i(t)\right\|},\;i=1,2$, are the two normalized deviation vectors at time $t$. SALI approaches zero exponentially fast in time (as a function of the largest or 2 largest Lyapunov exponents) for chaotic trajectories and non-zero, positive, values for regular trajectories \cite{skokos2004detecting}.

\begin{figure}
	\begin{subfigure}[tbhp]{.5\textwidth}
		\centering
		\includegraphics[width=0.9\linewidth]{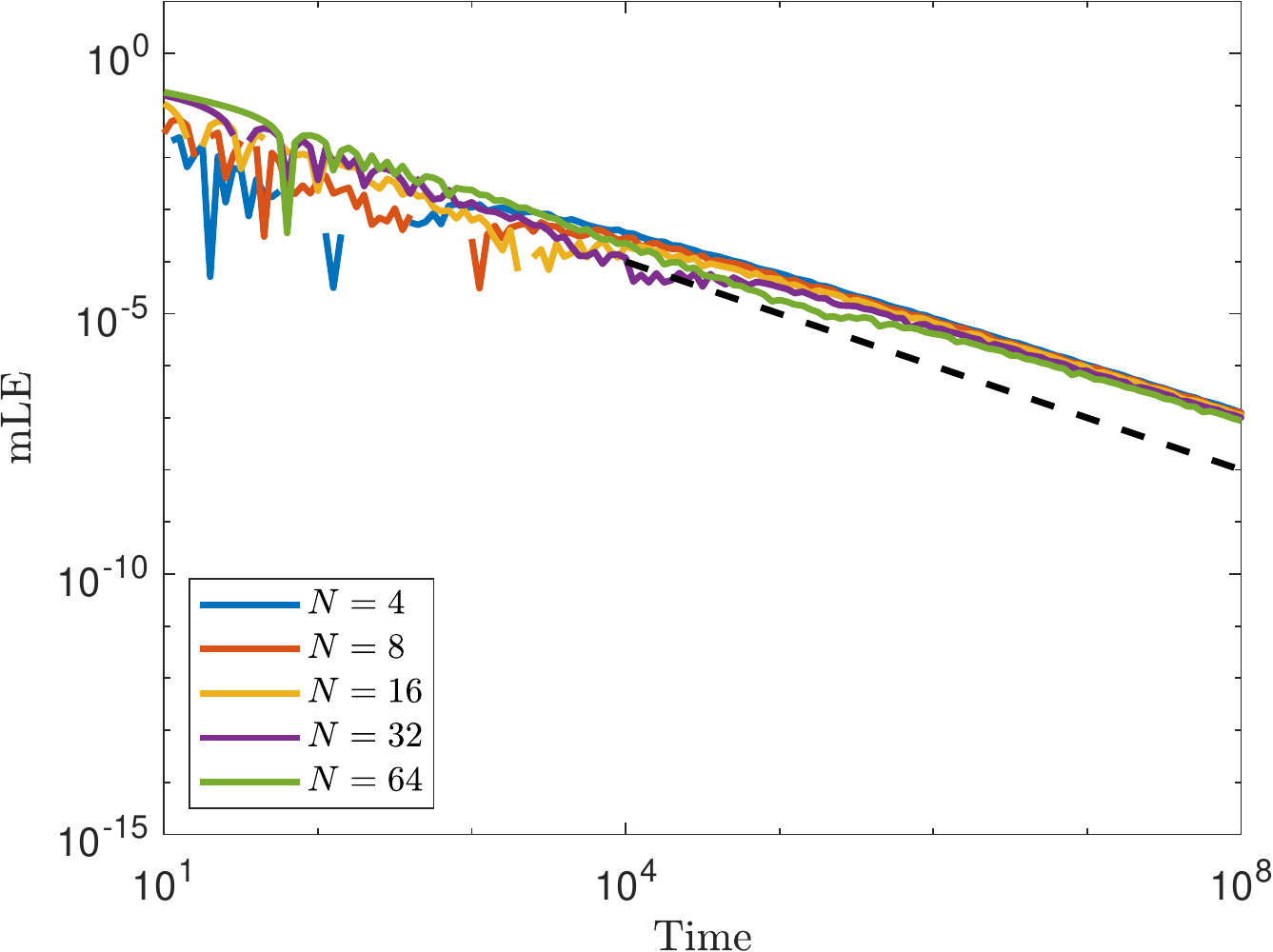}
		\caption{}
		\label{fig11:sfig1}
	\end{subfigure}%
	\begin{subfigure}[tbhp]{.5\textwidth}
		\centering
		\includegraphics[width=0.9\linewidth]{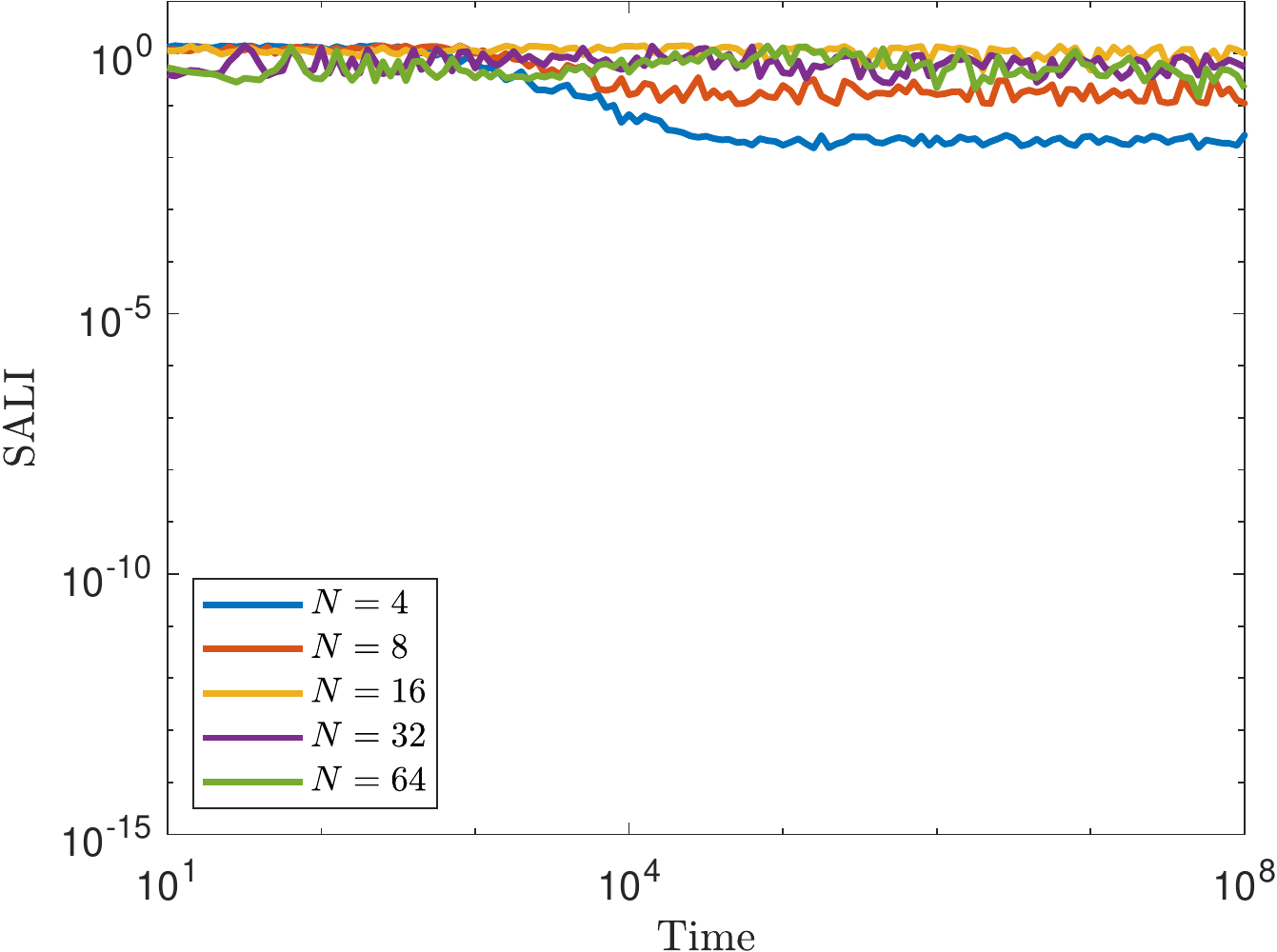}
		\caption{}
		\label{fig11:sfig2}
	\end{subfigure}
	\caption{Plot of mLE (panel a)) and SALI (panel b)) in time for a range of $N$ values seen in the insets (denoted by different colors) in the absence of variability, i.e., of the FPUT system \eqref{Hamilt}. Note that all axes are logarithmic. The black dashed line in panel (a) is the law $1/t$ of regular trajectories to guide the eye.}\label{fig:fig11}
\end{figure}

First, we consider the case without variability, that is the FPUT$-\alpha$ system \eqref{Hamilt}. We integrate the equations of motion \eqref{eqnmotions} and its corresponding variational equations (following Eq.\ \eqref{wode}) by using the tangent-map method \cite{skokos2010numerical} and Yoshida's fourth order symplectic integrator \cite{yoshida1990construction}. We have found that a time step of $0.01$ keeps the relative energy error below $10^{-9}$. In all our computations, the final integration time is $t=10^8$. Here, we use the same initial condition in Eq.\ \eqref{initxi} for all $N$. This initial condition then results in different energies for different $N$, i.e., $E=0.4775$ for $N=4$, $E=0.2714$ for $N=8$, $E=0.1447$ for $N=16$, $E=0.0747$ for $N=32$, and $E=0.0379$ for $N=64$. Our results in Fig.\ \ref{fig:fig11} show that all trajectories for $N=4,8,16,32,64$ are regular up to $t=10^8$, corroborated by the tendency of the mLEs to converge to zero following the $1/t$ law and SALI to tend to fixed positive values, shown in panels (a) and (b), respectively. These results are in agreement with the fact that energy recurrences in the homogeneous FPUT lattice \eqref{Hamilt} arise when it remains in the quasi-stationary state for extremely long times, making the approach to equipartition of energy unobservable.

\begin{figure}[htp]
	\centering
	\begin{subfigure}{.5\textwidth}
		\centering
		\includegraphics[width=0.9\linewidth]{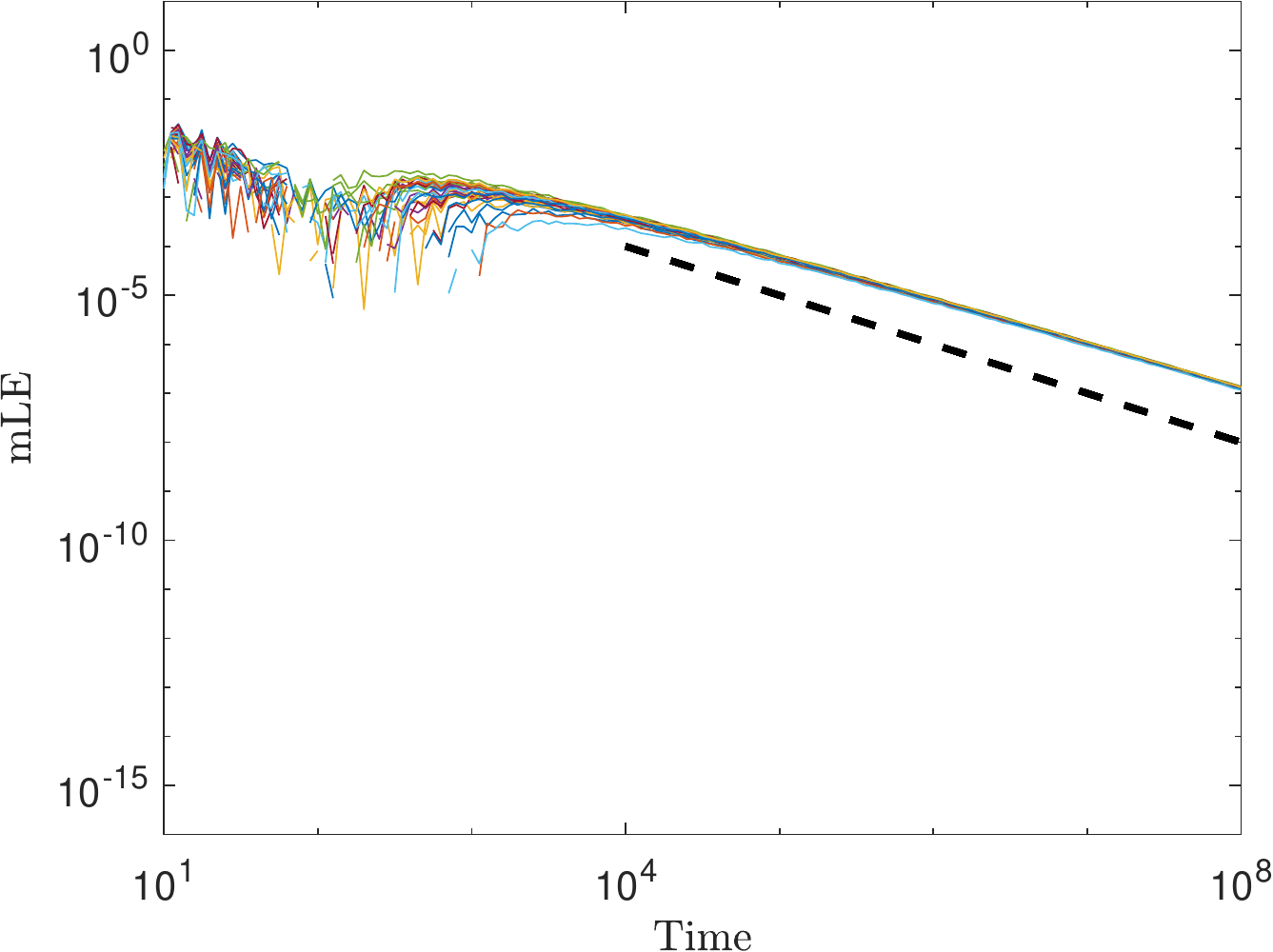}
		\caption{}
		\label{}
	\end{subfigure}%
	\begin{subfigure}{.5\textwidth}
		\centering
		\includegraphics[width=0.9\linewidth]{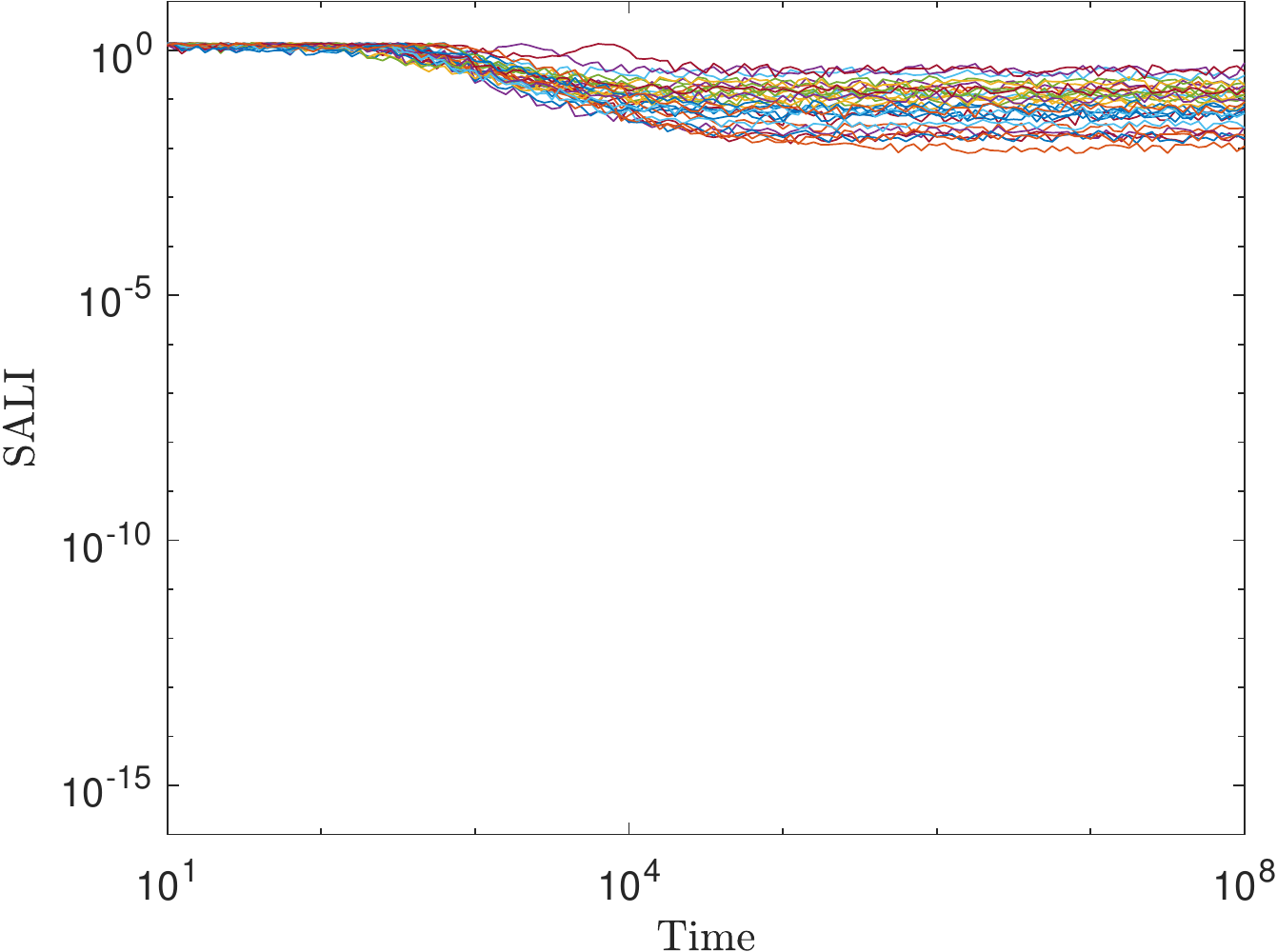}
		\caption{}
		\label{}
	\end{subfigure}\\
	\begin{subfigure}{.5\textwidth}
		\centering
		\includegraphics[width=0.9\linewidth]{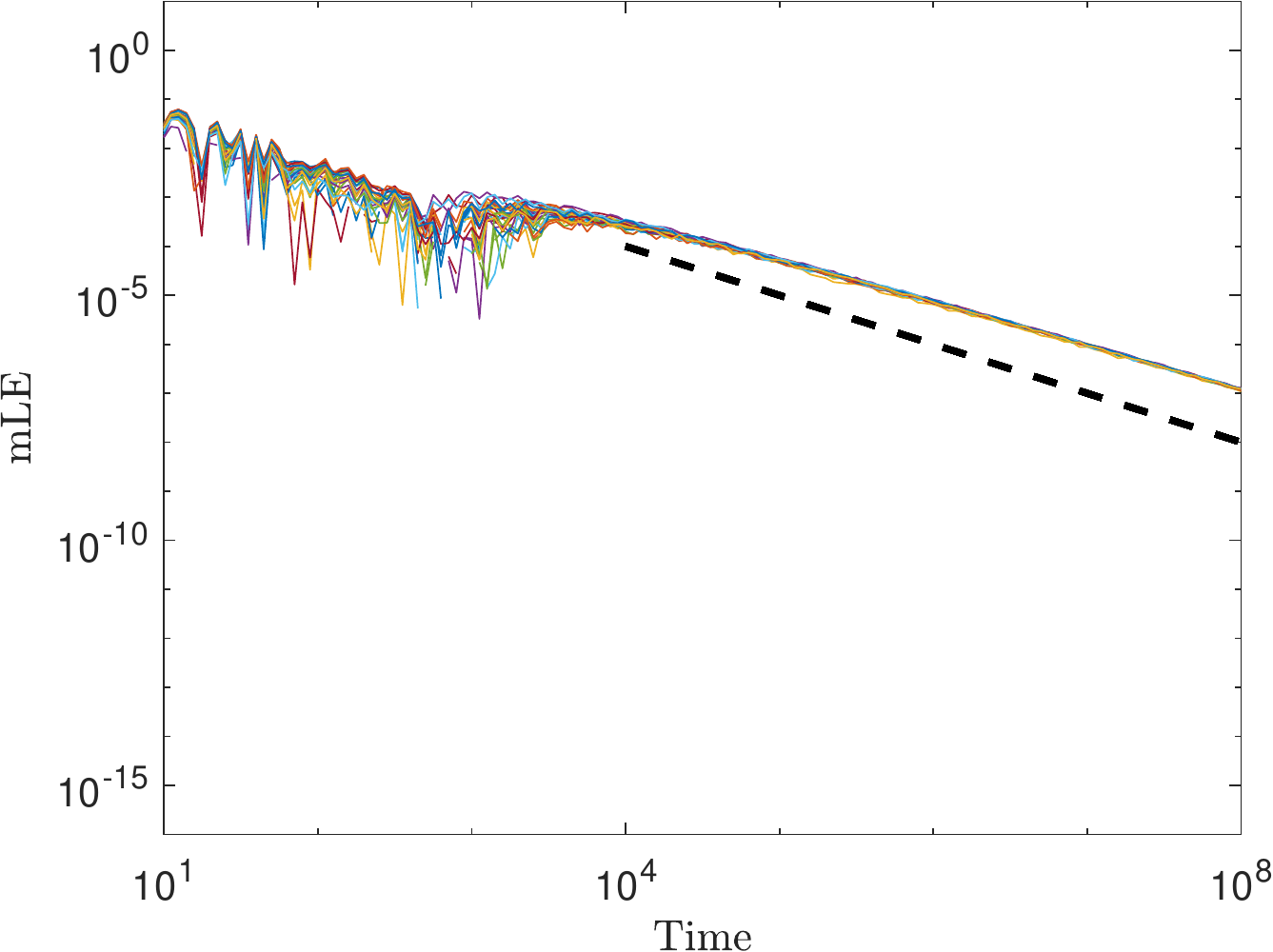}
		\caption{}
		\label{}
	\end{subfigure}%
	\begin{subfigure}{.5\textwidth}
		\centering
		\includegraphics[width=0.9\linewidth]{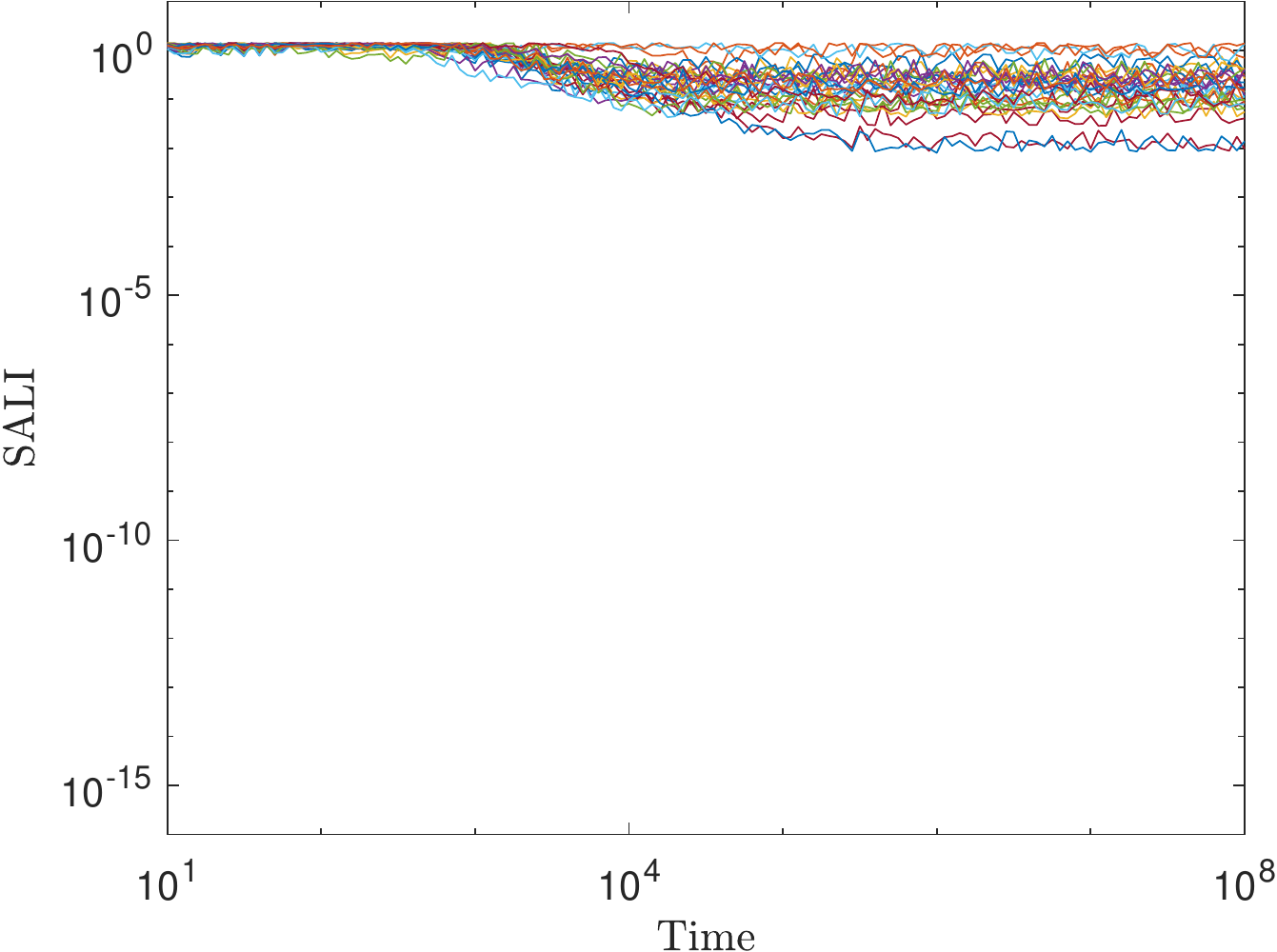}
		\caption{}
		\label{}
	\end{subfigure}\\
	\begin{subfigure}{.5\textwidth}
		\centering
		\includegraphics[width=0.9\linewidth]{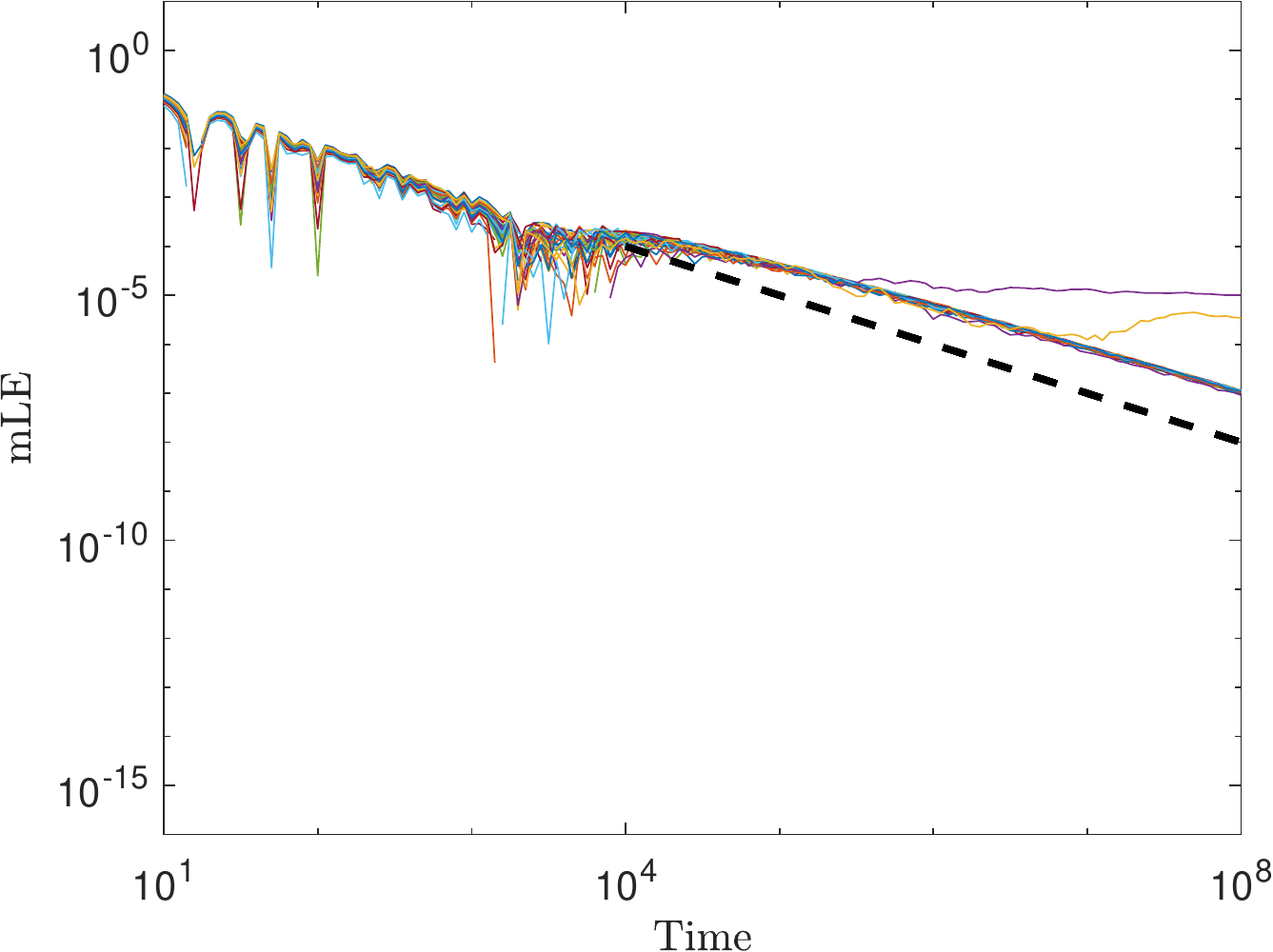}
		\caption{}
		\label{}
	\end{subfigure}%
	\begin{subfigure}{.5\textwidth}
		\centering
		\includegraphics[width=0.9\linewidth]{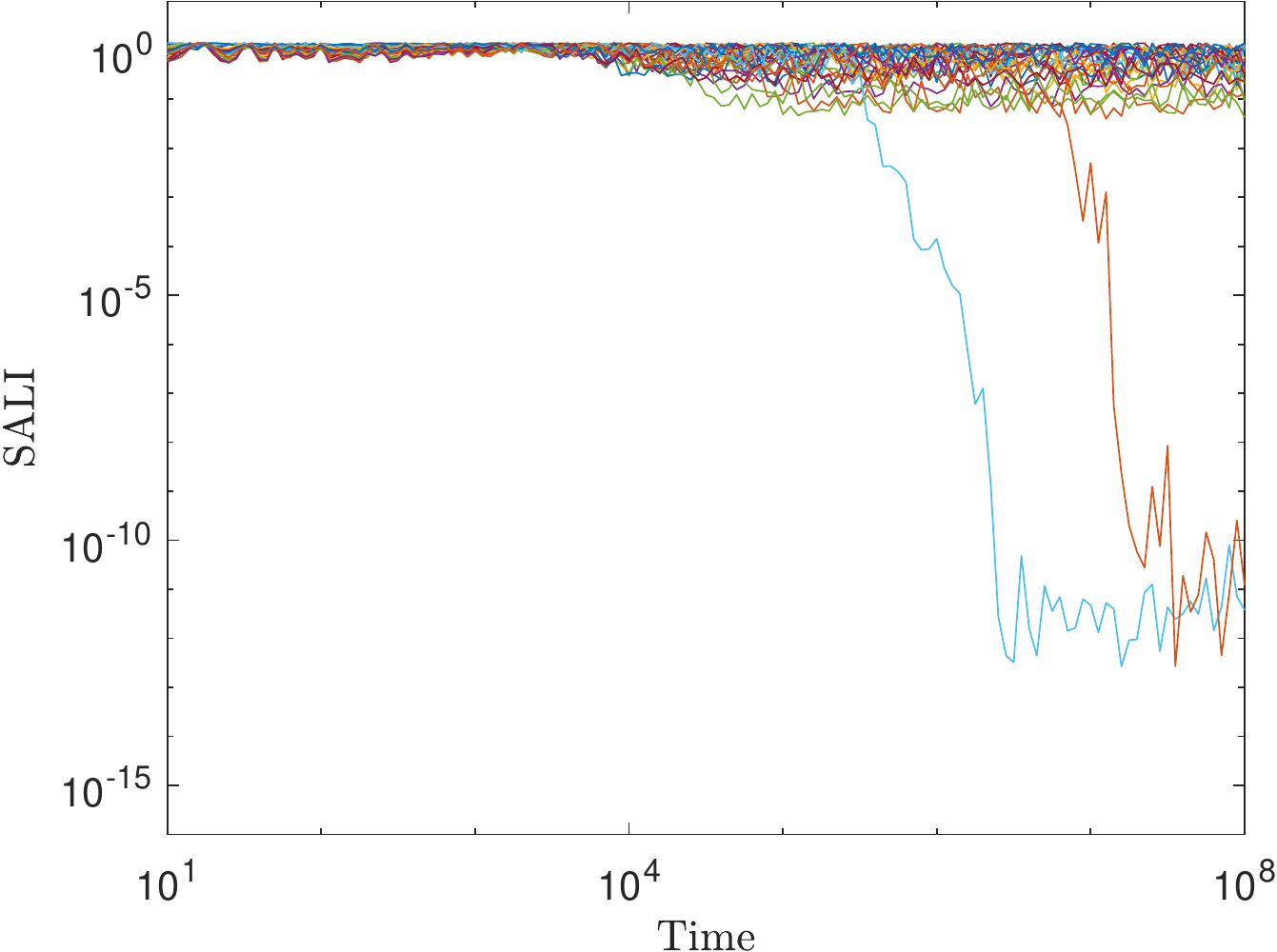}
		\caption{}
		\label{}
	\end{subfigure}
	\caption{Plot of mLE (panels (a), (c), (e)) and SALI (panels (b), (d), (f)) in time for 30 trajectories (denoted by different colors) and $\tau=10\%$ (see Eq.\ \eqref{varnonlin}). Panels (a), (b) are for $N=4$, panels (c), (d) for $N=8$ and panels (e), (f) for $N=16$. Note that all axes are logarithmic. The black dashed lines in panels (a), (c), (e) are the law $1/t$ of regular trajectories to guide the eye.}
	\label{fig:fig12}
\end{figure}

\begin{figure}[htp]
	\centering
	\begin{subfigure}{.5\textwidth}
		\centering
		\includegraphics[width=0.9\linewidth]{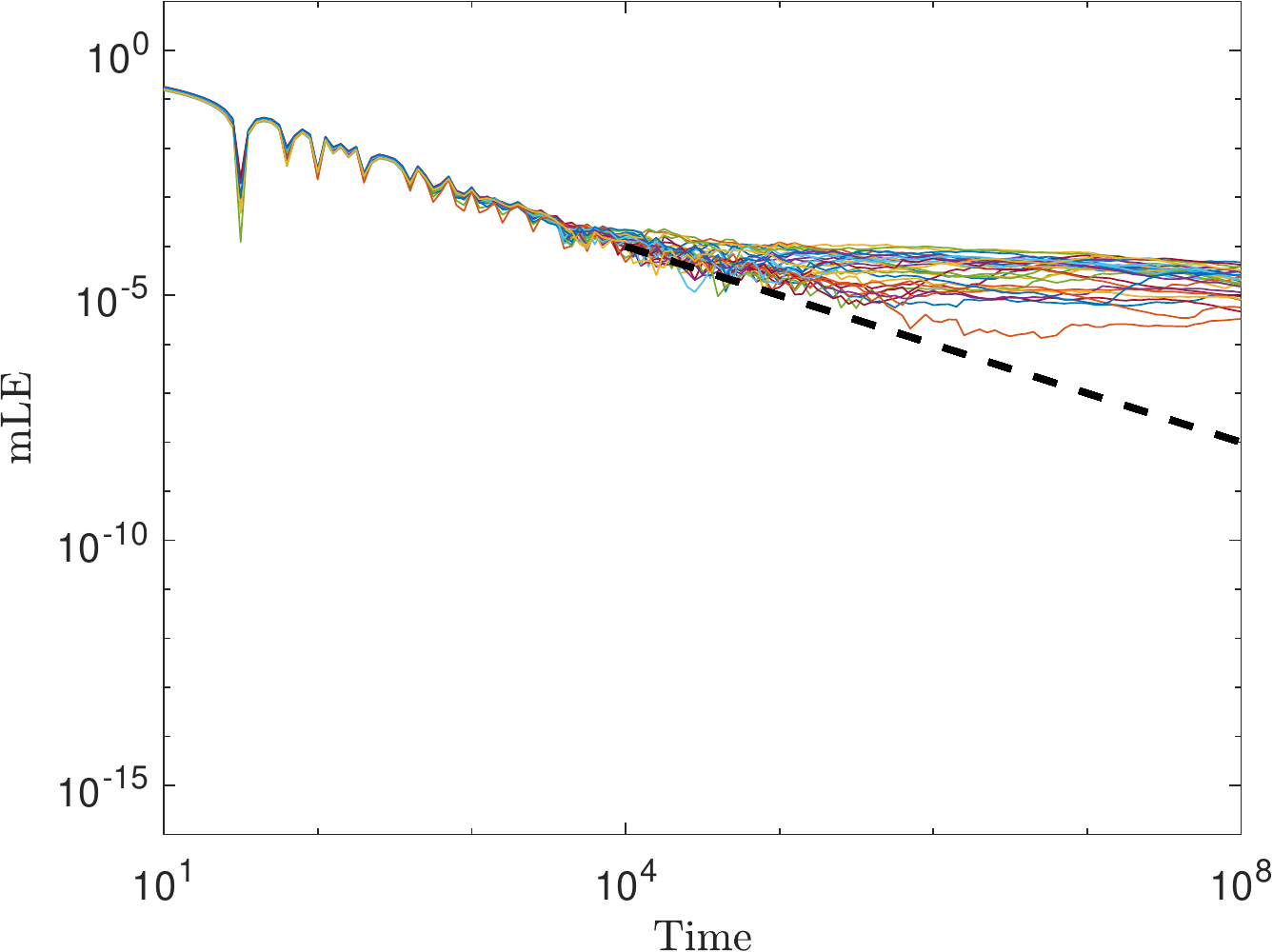}
		\caption{}
		\label{}
	\end{subfigure}%
	\begin{subfigure}{.5\textwidth}
		\centering
		\includegraphics[width=0.9\linewidth]{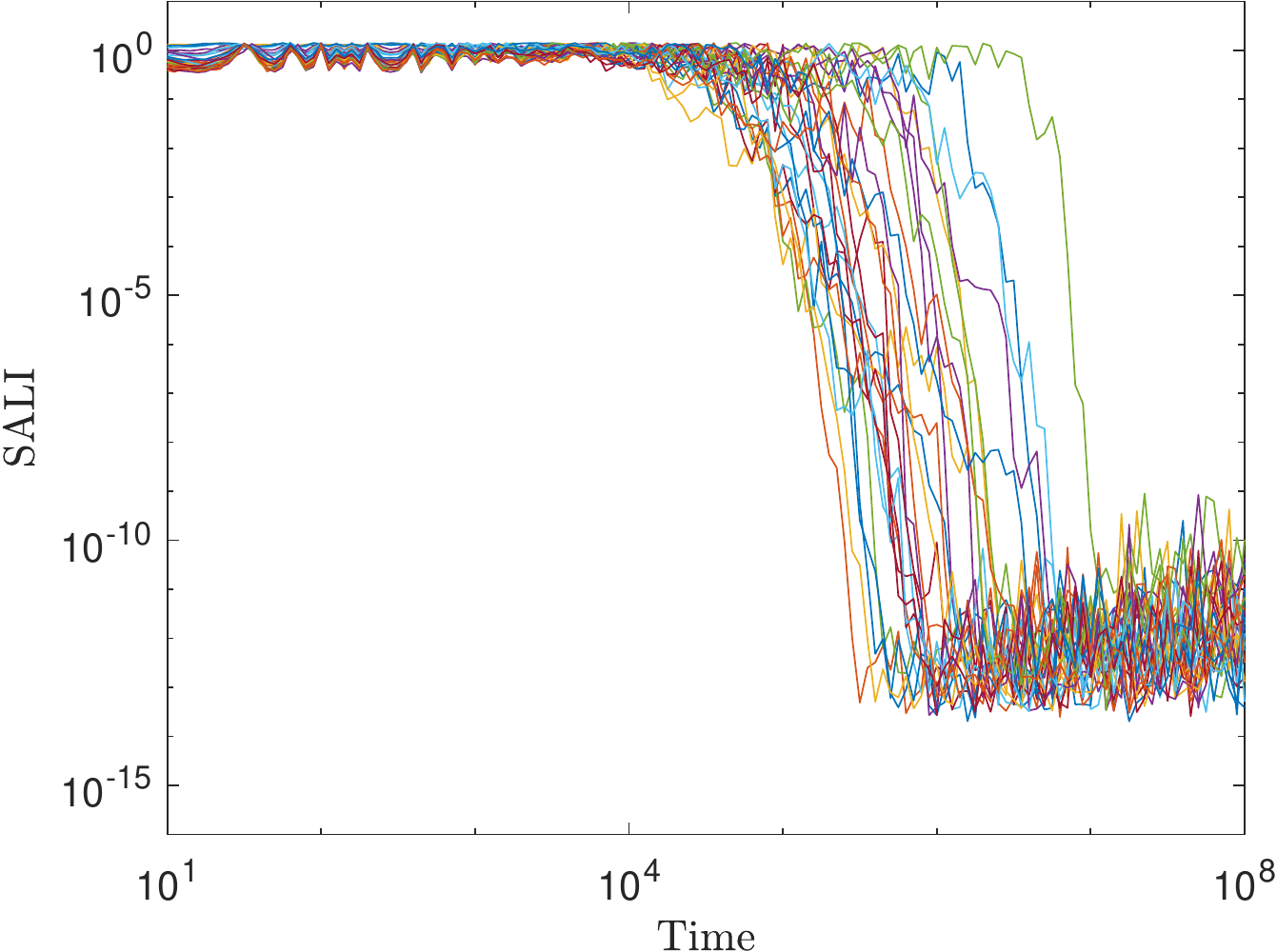}
		\caption{}
		\label{}
	\end{subfigure}\\
	\begin{subfigure}{.5\textwidth}
		\centering
		\includegraphics[width=0.9\linewidth]{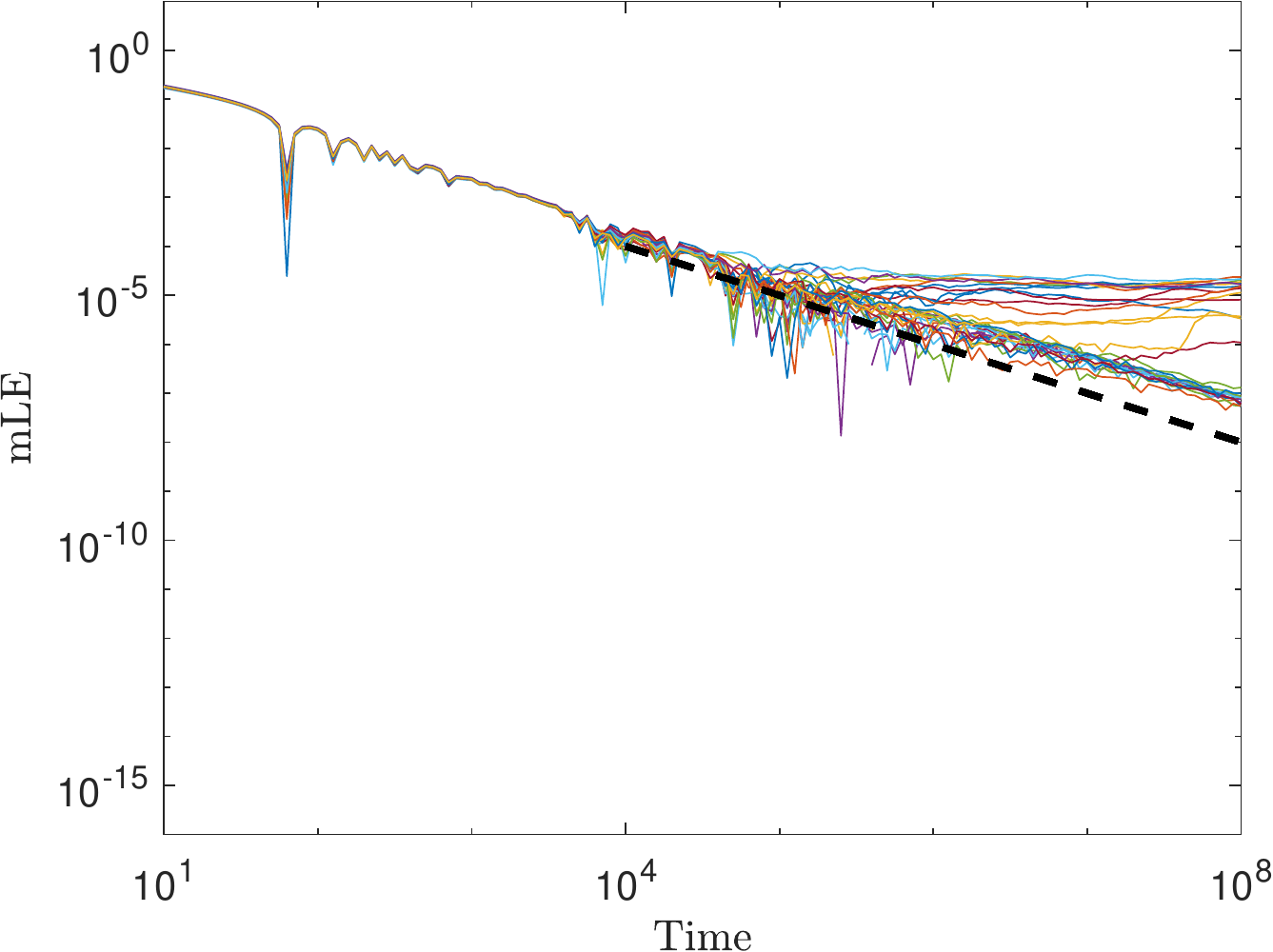}
		\caption{}
		\label{}
	\end{subfigure}%
	\begin{subfigure}{.5\textwidth}
		\centering
		\includegraphics[width=0.9\linewidth]{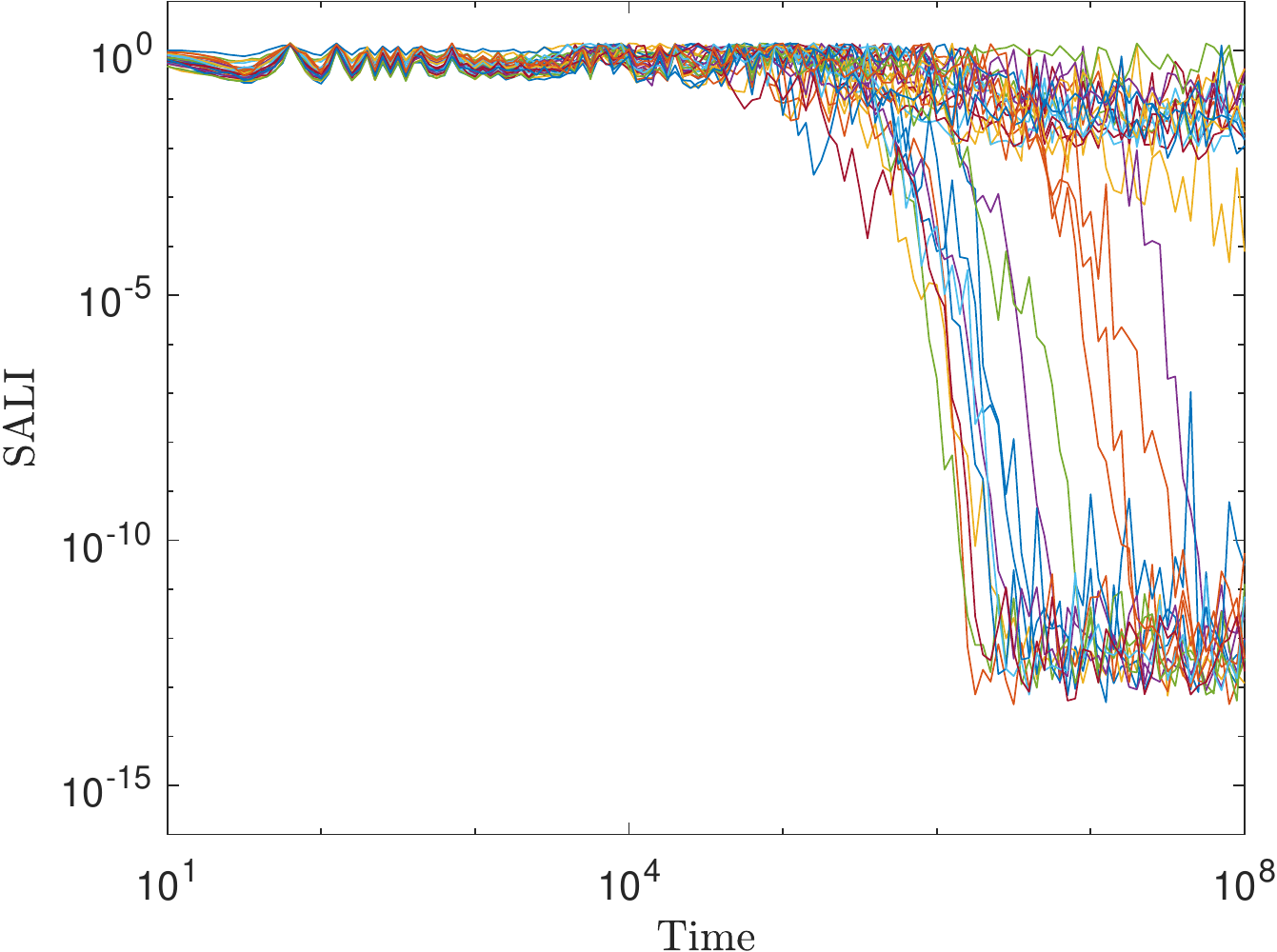}
		\caption{}
		\label{}
	\end{subfigure}
	\caption{Plot of mLE (panels (a), (c)) and SALI (panels (b), (d)) in time for 30 trajectories (denoted by different colors) and $\tau=10\%$ (see Eq.\ \eqref{varnonlin}). Panels (a), (b) are for $N=32$ and panels (c), (d) for $N=64$. Note that all axes are logarithmic. The black dashed lines in panels (a), (c) are the law $1/t$ of regular trajectories to guide the eye.}
	\label{fig:fig13}
\end{figure}

\begin{figure}[htp]
	\centering
	\includegraphics[width=0.5\textwidth]{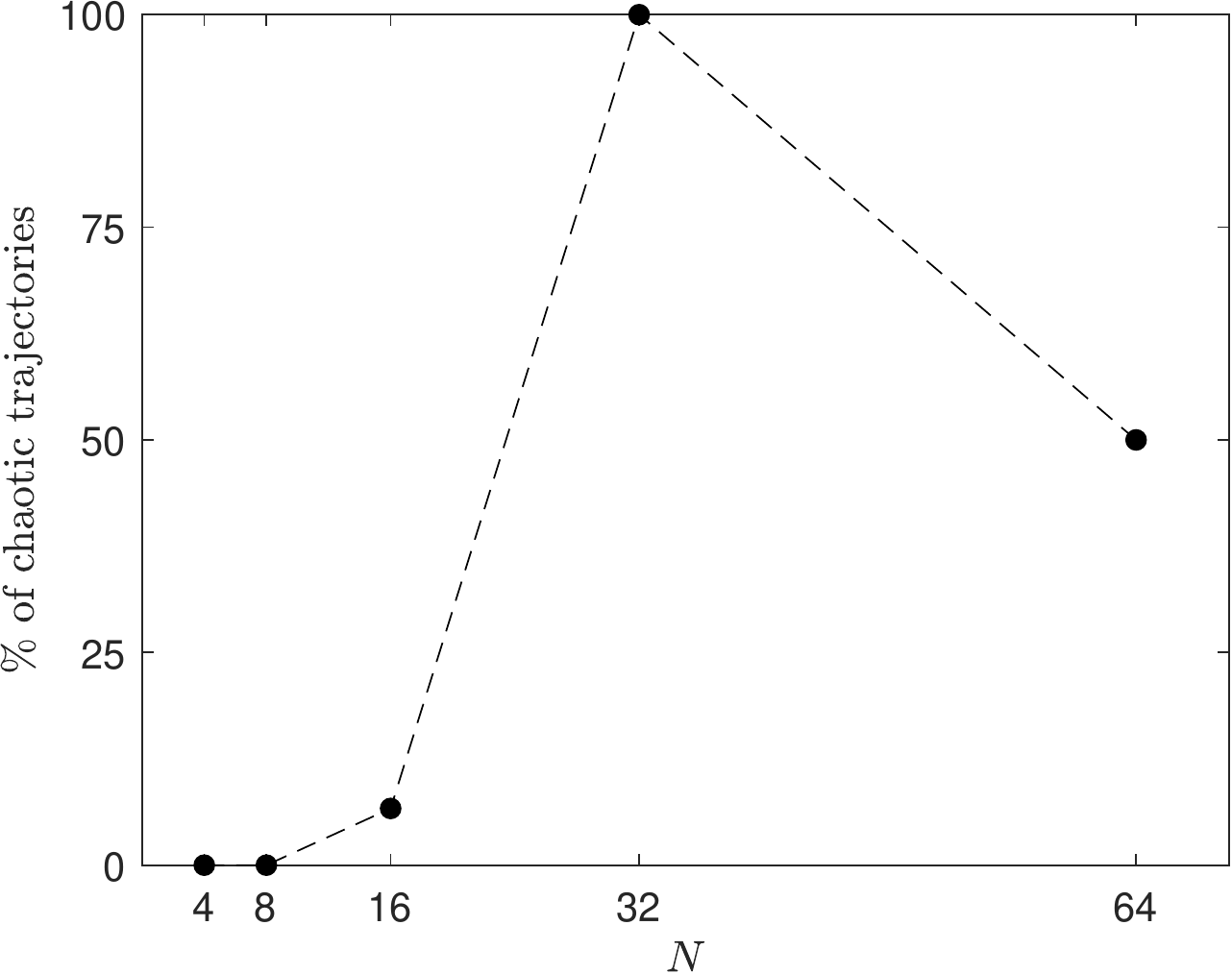}
	\caption{Percentage of chaotic trajectories as a function of $N$ for $30$ realisations of variability with $\tau=10\%$. The black-dash line segments connect the black points and are there to guide the eye.}\label{fig:fig14}
\end{figure}

Finally, we look at the case of $\tau=10\%<\tau_c$, for which we have observed almost energy localization in the first normal mode in Sec.\ \ref{sec3}. Since in this case we only know the equations of motion \eqref{varnonlin}, we integrated them using the DOP853 integrator \cite{bworld}, an explicit Runge-Kutta method of order $8$ due to Dormand and Prince, to achieve good numerical accuracy. We compute the chaotic indicators for $30$ realisations of the same percentage of variability $\tau=10\%$, while keeping the initial conditions fixed for each number of particles $N$. For $N=4$ and 8, all trajectories in panels (a)-(d) in Fig.\ \ref{fig:fig12} appear to be regular up to final integration time $t=10^8$, corroborated by the tendency of the mLEs to converge to zero following the $1/t$ law and SALI to tend to fixed positive values. However, for $N=16$, two of the 30 trajectories in panels (e), (f) in Fig.\ \ref{fig:fig13} are chaotic as their mLEs converge to positive values at $t=10^8$ and their SALI decrease to zero exponentially fast. Figure \ref{fig:fig13} shows that there are more chaotic orbits than those for the smaller values of $N$ in Fig.\ \ref{fig:fig12}. We show the percentage of chaotic trajectories (out of the 30 realisations) as a function of $N$ in Fig.\ \ref{fig:fig14}, where the increase from $N=4,8,16$ to $N=32,64$ is apparent. These results suggest that in the case of almost complete energy localization, variability promotes chaos in the system as the number of particles increases. However, further studies are required to determine whether the increase is monotone.

\section{Conclusions and discussion}\label{conc}

In this paper, we have considered a disordered FPUT$-\alpha$ system with variations in its parameters (also called variability) to take into account inherent manufacturing processes. By using a two normal-mode approximation, we have been able to explain the mechanism for energy localization and blow up of solutions for percentage of variability bigger than a threshold, that we have been able to compute using our theory. Moreover, we have also studied the effect of variability in the chaotic behavior of the  system calculating the maximum Lyapunov exponent and Smaller Alignment Index for a number of realizations for the same variability percentage that corresponds to energy-localization. We have found that, when there is almost energy localization, it is more frequent for the trajectories to be chaotic with the increase of the number of particles $N$ for the same percentage of variability, smaller than the threshold.

Finally, while it has been shown previously that variability leads to energy-recurrence breakdown and energy localization, we have also shown here that by increasing the percentage of variability beyond a threshold that we determined using our theory, the solutions of the system may blow up in finite-time. This is because we have started with the equations of motion without a Hamiltonian that would allow us to keep the energy of the system constant \cite{porter2018}. 
The case of the Hamiltonian model with heterogeneity, cf.\ Eq.\ \eqref{dis1}, will be studied in a future publication. 

\section*{Credit authorship contribution statement}
Zulkarnain: Investigation, Writing -- Original Draft. H.\ Susanto: Conceptualization, Supervision, Methodology, Writing -- review \& editing. C.\ Antonopoulos: Conceptualization, Supervision, Methodology, Writing -- review \& editing.

\section*{Declaration of Competing Interest}
The authors declare that they have no known competing financial interests or personal relationships that could have appeared to influence the work reported in this paper.

\section*{Acknowledgement} Z is supported by the Ministry of Education, Culture, Research, and Technology of Indonesia through a PhD scholarship (BPPLN). HS is supported by Khalifa University through a Faculty Start-Up Grant (No.\ 8474000351/FSU-2021-011) and a Competitive Internal Research Awards Grant (No.\ 8474000413/CIRA-2021-065). The authors acknowledge the use of the High Performance Computing Facility (Ceres) and its associated support services at the University of Essex in the completion of this work. The authors are also grateful to the referees for their comments and feedback that improved the manuscript. 
\newpage
%
%
%

\end{document}